\documentclass[a4paper,10pt]{article}

\usepackage{amsmath,amssymb}
\usepackage{amsthm}
\usepackage{palatino}
\usepackage{amscd}
\usepackage{graphicx}
\usepackage{epsfig}
\usepackage{verbatim}
\usepackage{calc}

\title{\textbf{Parameters for which the Lawrence-Krammer representation is reducible}}

\author{Claire Levaillant and David Wales\\
cl@caltech.edu, dbw@caltech.edu\\\\
Caltech, MC 253-37, Pasadena, CA 91125}

\newcommand{\ot}{\otimes_{_{H}}}

\newcommand{\m}{\mu}
\newcommand{\n}{\nu}
\newcommand{\Q}{\mathbb{Q}}
\newcommand{\R}{\mathbb{R}}
\newcommand{\be}{\beta}
\newcommand{\al}{\alpha}
\newcommand{\xb}{x_{\beta}}
\newcommand{\e}{\epsilon}
\newcommand{\lra}{\longrightarrow}

\newcommand{\la}{\lambda}

\newcommand{\unsurr}{\frac{1}{r}}
\newcommand{\unsur}{\frac{1}}

\newcommand{\ali}{\al_i}
\newcommand{\xali}{x_{\al_i}}

\newcommand{\U}{\mathcal{U}}
\newcommand{\V}{\mathcal{V}}
\newcommand{\IH}{\mathcal{H}}

\newcommand{\lb}{\lbrace}
\newcommand{\rb}{\rbrace}
\newcommand{\noin}{\noindent}
\newcommand{\bigp}{\Big(}
\newcommand{\bigpd}{\Big)}

\newcommand{\ih}{\IH_{F,r^2}}
\newcommand{\D}{\mathcal{D}}

\newcommand{\W}{\mathcal{W}}
\newcommand{\di}{\text{dim}}

\newcommand{\da}{\downarrow}

\newcommand{\B}{\mathcal{B}}

\newcommand{\tr}{\triangle}

\newcommand{\eg}{&=&}
\newcommand{\g}{\gamma}

\newcommand{\cil}{\frac{n(n-3)}{2}}
\newcommand{\dbw}{\frac{(n-1)(n-2)}{2}}
\newcommand{\chl}{\frac{n(n-1)}{2}}

\newcommand{\mfl}{\unsur{r^{2n-3}}}
\newcommand{\jca}{\unsur{r^{n-3}}}

\begin{document}
\maketitle \hrule
\section*{Abstract}
We show that the representation, introduced by Lawrence and Krammer
to show the linearity of Braid groups, is generically irreducible,
but for that for some values of its two parameters when these are
specialized to complex numbers, it becomes reducible. To do so, we
construct a representation of degree $\chl$ of the BMW algebra of
type $A_{n-1}$ inside the Lawrence-Krammer space. As a
representation of the Braid group on $n$ strands, it is equivalent
to the Lawrence-Krammer representation where the two parameters of
the BMW algebra are related to the two parameters of the
Lawrence-Krammer representation. We give the values of the
parameters for which the representation is reducible and give the
proper invariant subspaces in some cases. We use this representation
to show that for these special values of the parameters and other values, the BMW algebra of type $A_{n-1}$ is not semisimple.\\
\hrule
\section{Introduction}
\subsection{Introduction and main results}
In \cite{KR}, Daan Krammer constructed a representation of the Braid
group in order to show that it is linear. Since this representation
was earlier introduced by Ruth Lawrence in \cite{RL}, it is called
the Lawrence-Krammer representation. In this paper, we examine a
representation of degree $\chl$ of the BMW algebra of type $A_{n-1}$
in the Lawrence-Krammer space. As a representation of the Braid
group on $n$ strands, it is equivalent to the Lawrence-Krammer
representation (abbreviated L-K representation). By studying this
representation we show that the L-K representation is generically
irreducible but that for some values of its two parameters when
these are specialized to complex numbers, it becomes reducible.
Throughout the paper, we let $l$, $m$ and $r$ be three nonzero
complex parameters, where $m$ and $r$ are related by $m=\unsurr-r$.
We define $\ih(n)$ as the Iwahori-Hecke algebra of the symmetric
group $Sym(n)$ over the field $F=\Q(l,r)$ with generators
$g_1,\dots,g_{n-1}$, that satisfy the Braid relations and the
relation $g_i^2+m\,g_i=1$ for all $i$. Our definition is the same as
the definition of \cite{M} after the generators have been rescaled
by a factor $\unsurr$. Our main result is as follows.

\newtheorem{theo}{Theorem}
\begin{theo}(Main theorem)\hfill\\
Let $n$ be an integer with $n\geq 3$ and let $m$, $l$ and $r$ be
three nonzero complex parameters, where $m$ and $r$ are related by
$m=\unsurr-r$. Assume that $\ih(n)$ is semisimple, and so assume
that $r^{2k}\neq 1$ for every integer $k\in\lb 1,\dots,n\rb$.
\\When $n\geq 4$, the Lawrence-Krammer representation of the BMW
algebra of type $A_{n-1}$ with parameters $l$ and $m$ over the field
$\Q(l,r)$ is irreducible, except when $l\in\lb r,-r^3,
\unsur{r^{2n-3}},\jca, -\jca\rb$, when it is reducible.
\\When $n=3$, the Lawrence-Krammer representation of the BMW algebra
of type $A_2$ with parameters $l$ and $m$ over the field $\Q(l,r)$
is irreducible, except when $l\in\lb -r^3,\unsur{r^3},1,-1\rb$, when
it is reducible.
\end{theo}

A consequence of this result and of the method that we use is the
following.

\begin{theo}\hfill\\
Let $n$ be an integer with $n\geq 4$ and let $l$, $m$
and $r$ be three nonzero complex parameters, where $m$ and $r$ are related by $m=\unsurr-r$.\\
Suppose $n\geq 4$. If $r^{2k}=1$ for some $k\in\lb 2,\dots,n\rb$ or
if $l$ belongs to the set of values $\lb
r,-r^3,\jca,-\jca,\mfl,-r^{2n-3},r^{n-3},-r^{n-3},\unsur{r^3},-\unsurr\rb$,
the BMW algebra of type $A_{n-1}$ with parameters $l$ and $m$ over
the field $\Q(l,r)$ is not semisimple.\\
(Case $n=3$). If $r^4=1$ or $r^6=1$ or if $l\in\lb -r^3,
\unsur{r^3},1,-1\rb$, the BMW algebra of type $A_2$ with parameters
$l$ and $m$ over the field $\Q(l,r)$ is not semisimple.
\end{theo}

\noin In \cite{W}, Wenzl states that the BMW algebra of type
$A_{n-1}$ is semisimple except possibly if $r$ is a root of unity or
$l=r^n$, for some $n\in\mathbb{Z}$. Here, Theorem $2$ gives
instances of when the algebra is not semisimple. The result of this
theorem is also contained in the recent work of Hebing Rui and Mei
Si (see \cite{HR}). They use the representation theory of cellular
algebras.

\subsection{The method}
We show that the action on a proper invariant subspace of the
Lawrence-Krammer space must be an Iwahori-Hecke algebra action.

First, we study the Iwahori-Hecke algebra representations of small
degrees and investigate wether they may occur inside the L-K space
and if so for which values of $l$ and $r$. We will denote by
$\V^{(n)}$ the L-K space. We show that if there exists a
one-dimensional invariant subspace inside $\V^{(n)}$, it forces the
value $\unsur{r^{2n-3}}$ for $l$, except when $n=3$ when it forces
$l\in\lb -r^3,\unsur{r^3}\rb$. Conversely, for these values of $l$
and $r$, there exists a one-dimensional invariant subspace of
$\V^{(n)}$ and the representation is thus reducible. Similarly, we
show that if there exists an irreducible $(n-1)$-dimensional
invariant subspace inside $\V^{(n)}$, it forces $l=\unsur{r^{n-3}}$
or $l=-\unsur{r^{n-3}}$ in the case when $n\neq 4$ and $l\in\lb
-r^3,\unsurr,-\unsurr\rb$ in the case when $n=4$. Conversely, for
each of these values of $l$ and $r$, there exists an irreducible
$(n-1)$-dimensional subspace of $\V^{(n)}$, which shows the
reducibility of the representation in these cases as well.

Second, we identify a proper invariant subspace of $\V^{(n)}$ which
is nontrivial when $l=r$ (case $n\geq 4$) or $l=-r^3$. This shows
that the representation is also reducible in these cases.

Third, we study in detail the small cases $n\in\lb 3,4,5,6\rb$.

At last, when $n\geq 7$, we use a result from representation theory:
the irreducible representations of $\ih(n)$ have degrees $1, n-1,
\cil,\dbw$ or degree greater than $\dbw$, except in the case $n=8$,
when they have degrees $1,7,14,20,21$ or degrees greater than $21$.
We use this fact, and proceed by induction on $n\geq 5$ to show that
if $\V^{(n)}$ is reducible, it forces $l\in\lb
r,-r^3,\unsur{r^{2n-3}},\jca,-\jca\rb$. To do so, we use the fact
that if the dimension of a proper invariant subspace $\W$ of
$\V^{(n)}$ is large enough, then the intersections
$\W\cap\V^{(n-1)}$ and $\W\cap\V^{(n-2)}$ are nontrivial.
\subsection{Definitions}
\subsubsection{The BMW algebra}
We recall below the defining relations of the BMW algebra
$B(A_{n-1})$ (or simply $B$) of type $A_{n-1}$ with nonzero complex
parameters $l$ and $m$ over the field $\Q(l,r)$, where $r$ is a root
of the quadratic $X^2-mX+1$. This algebra has two sets of $(n-1)$
elements, namely the invertible $g_i$'s that satisfy the Braid
relations $(1)$ and $(2)$ and generate the algebra and the $e_i$'s
that generate an ideal. For nodes $i$ and $j$ with $1\leq i,j\leq
n-1$, we will write $i\sim j$ if $|i-j|=1$ and $i\not\sim j$ if
$|i-j|>2$.
\begin{eqnarray}
g_ig_j & = & g_jg_i \;\;\;\;\;\;\;\;\;\;\;\;\;\;\;\;\;\;\;\;\;\;\;\;\;\;\;\,\;\;\;\;\; \text{if $i\not\sim j$}\\
g_ig_jg_i &=& g_jg_ig_j \;\;\;\;\;\;\;\;\;\;\;\;\;\;\;\;\;\;\;\;\;\;\;\;\;\;\;\;\; \text{if $i\sim j$}\\
e_i & = & \frac{l}{m}*(g_i^2+m\;g_i-1)\qquad\text{for all $i$}\\
g_ie_i &=& l^{-1}e_i \qquad\qquad\qquad\qquad\;\;\text{for all $i$}\\
e_ig_je_i & = & le_i
\;\;\;\;\;\;\;\;\;\;\;\;\;\;\;\;\;\;\;\;\;\;\;\;\;\;\;\;\;\;\;\;\;\;\;\!\text{if
$i\sim j$}
\end{eqnarray}
We will also use some direct consequences of these defining
relations (see \cite{CGW}, Proposition $2.1$):
\begin{eqnarray}
e_ig_i & = & l^{-1}e_i\qquad\qquad\qquad\qquad\qquad\!\!\!\text{for all $i$}\\
g_i^2 & = & 1-mg_i+ml^{-1}e_i\qquad\qquad\;\,\!\!\!\text{for all $i$}\\
g_i^{-1} & = & g_i + m-m\;e_i\qquad\qquad\qquad\!\!\!\text{for all
$i$}
\end{eqnarray}
as well as the following "mixed Braid relations" (see \cite{CGW},
Proposition $2.3$):
\begin{eqnarray}
g_ig_je_i & = & e_je_i\;\;\;\;\;\;\;\;\;\;\;\;\;\;\;\;\;\;\;\;\;\;\;\;\;\;\;\;\;\;\;\text{if $i\sim j$}\\
g_ie_je_i & = & g_je_i + m(e_i-e_je_i)\;\;\;\;\;\;\;\,\text{if
$i\sim j$}
\end{eqnarray}
This algebra was shown by Morton and Wassermann to be isomorphic to
the tangle algebra of Morton and Traczyk (see \cite{MW}). All the
algebraic relations given in this paper have a geometric formulation
in terms of tangles. In particular, we will use the tangles in
$\S3.4$.
\subsubsection{The Lawrence-Krammer space}
We now recall some terminology associated with root systems of type
$A_{n-1}$. Let $M=(m_{ij})_{1\leq i\leq j\leq n-1}$ be the
Coxeter matrix of type $A_{n-1}$.\\
Let $(\al_1,\dots,\al_{n-1})$ be the canonical basis of $\R^{n-1}$
and let's define a bilinear form $B_M$ over $\R^{n-1}$ by:
$$B_M(\al_i,\al_j)=-cos\bigg(\frac{\pi}{m_{ij}}\bigg)$$
By the theory in \cite{BOU}, $B_M$ is an inner product that we will
simply denote by $(\,|\,)$. Let $r_i$ denote the reflection with
respect to the hyperplane Ker$(\al_i|.)$ of $\R^{n-1}$ and so:
$$\forall x\in\R^{n-1},\,r_i(x)=x-2(\al_i|x)\,\al_i$$
Finally, let $\phi^{+}$ denote the set of $\chl$ positive roots:
\begin{equation*}\begin{split}\phi^{+}=\lb \al_1,\al_2,\al_2+\al_1,&\al_3,\al_3+\al_2,\al_3+\al_2+\al_1,\dots,\\
&\al_{n-1},\al_{n-1}+\al_{n-2},\al_{n-1}+\al_{n-2}+\dots+\al_1\rb\end{split}\end{equation*}
We define $\V^{(n)}$ as the vector space over $\Q(l,r)$ with basis
the $x_{\be}$'s, indexed by the positive roots $\be\in\phi^{+}$.
Thus, $\di_F\V^{(n)}=|\phi^{+}|=\chl$.
The so-defined space $\V^{(n)}$ is the Lawrence-Krammer space.\\
To each positive root, we associate an element of the BMW algebra in
the following way:
\begin{itemize}\item To $\al_1$ we associate $e_1$.
\item To $\al_i=r_{i-1}\dots\,r_1\,r_i\dots\,r_2(\al_1)$, we
associate the algebra element\\ $g_{i-1}\dots g_1g_i\dots g_2e_1$,
which after using the defining rules $(1)$ and $(9)$ above
simplifies to $e_i\dots e_1$.
\item To $\al_j+\dots+\al_i=r_j\dots r_{i+1}(\al_i)$ where $j\geq
i+1$, we associate the algebra element $g_j\dots g_{i+1} e_i\dots
e_1$.
\end{itemize}

\section{The representation}
\subsection{The BMW left module}
In what follows, $F$ still denotes the field $\Q(l,r)$ and $H$
denotes the Hecke algebra of the symmetric group $Sym(n-2)$ over the
field $F$ with generators $g_3,\dots,g_{n-1}$ and relations the
Braid relations and the relations $g_i^2+m\,g_i=1$ for each $i$. As
$r^2+m\,r-1=0$, our base field $F$ is a one-dimensional $H$-module
for the action given by $g_i.1=r$ for every integer $i$ with $3\leq
i\leq n-1$. We define $B_1$ as the quotient of two left ideals of
$B$:
$$B_1=Be_1/<Be_ie_1>_{i=3\dots n-1}$$
Since $e_i$ commutes with $e_j$ for any $i\not\sim j$, we have for
each node $i$ with $3\leq i\leq n-1$:
$$e_1(g_i^2+m\,g_i-1)=\frac{l}{m}e_1e_i=0\;\;\text{in}\;B_1$$ Then, $B_1$ is a right
$H$-module. Thus, $B_1$ is a left $B$-module and a right $H$-module.
Since $F$ is an $H$-module, we get a left $B$-module by considering
the tensor product
$$B_1\otimes_{H}F$$
This $B$-module is precisely the left representation of $B$ that we
study in this paper. Its degree is $\chl$ since by the forthcoming
computations, we have:
\begin{equation*}\begin{split}B_1=Span_F(e_1,e_2e_1,&g_2e_1,e_3e_2e_1,g_3e_2e_1,g_3g_2e_1,\dots,\\
&e_{n-1}\dots e_1,g_{n-1}e_{n-2}\dots e_1,g_{n-1}\dots
g_2e_1).H\end{split}\end{equation*} We will denote by $\B$ the
spanning set above. The action of the $g_i$'s on the elementary
tensors $b\ot 1$, where $b\in\B$, was computed in section $2.2$
below. These computations show in particular that if $G_1(n)$
denotes the matrix of the left action of $g_1$ on the vectors
\begin{multline*}e_1\ot 1,e_2e_1\ot 1,g_2e_1\ot 1,e_3e_2e_1\ot
1,g_3e_2e_1\ot 1,g_3g_2e_1\ot 1,\dots,\\e_{n-1}\dots e_1\ot
1,g_{n-1}e_{n-2}\dots e_1\ot 1,g_{n-1}\dots g_2e_1\ot
1\end{multline*} of $B_1\ot F$, and if $\text{det}\,G_1(n)$ denotes
its determinant, we have:
\begin{eqnarray*}
\text{det}\,G_1(3)&=&-\unsur{l}\\
\text{det}\,G_1(n)&=&-r^{n-3}\,\text{det}\,G_1(n-1)\qquad\forall
n\geq 4
\end{eqnarray*}
Thus, the determinant of $G_1(n)$ is nonzero, which shows that these
vectors are linearly independent.
\\
We notice that there is a bijection between $\B$ and the set of
positive roots $\phi^{+}$, as described in the previous section.
Let's name this bijection $u$.

\subsection{The action by the $g_k$'s}
We describe further the representation by computing the action of
the $g_k$'s on the elementary tensors $b\otimes_H 1$, where $b$ is
an algebra element in the spanning set $\B$. An element $b$ of $\B$
is of the form:\newcounter{ctr}\setcounter{ctr}{1}
$$g_j\dots g_{i+1}e_i\dots e_1\;\;\text{with}\;\;j>i\geq 1\;\;(\text{or simply}\;\;g_{j,i+1}e_{i,1})\qquad\qquad\qquad\;(\Roman{ctr})$$ that we will
refer to as of type \newcounter{c}\setcounter{c}{1}$(\Roman{c})$, or
of the form:\newcounter{compteur}\setcounter{compteur}{2}
$$e_i\dots e_1\;\;\text{with}\;\;i\geq 1\;\;(\text{or simply}\;\;e_{i,1})\qquad\qquad\qquad\,\,\qquad\qquad\qquad\qquad(\Roman{compteur})$$ referred to as of type
\setcounter{c}{2}$(\Roman{c})$. \\For $i\geq j$, we set
$g_{i,j}=g_i\dots g_j$ and $e_{i,j}=e_i\dots e_j$, where $g_{i,i}$
and $e_{i,i}$ are simply $g_i$ and $e_i$ respectively. When $i<j$,
we define $g_{i,j}$ to be the identity.\\
In what follows, we fix $i$ and $j$ as in \setcounter{c}{1}
$(\Roman{c})$ and \setcounter{c}{2} $(\Roman{c})$. There are several
cases.

\subsubsection{Action by $g_{i-1}$ (Case $A$)} Let's first compute the action of $g_{i-1}$ for $i\geq
2$ on elements of both types. We have for the type
\setcounter{c}{1}$(\Roman{c})$: \begin{eqnarray*}
g_{i-1}.b&=&g_{j,i+1}g_{i-1}e_{i,1}\qquad\qquad\qquad\qquad\qquad\qquad\;\;\,\text{by $(1)$}\\
&=& g_{j,i}e_{i-1,1} +
m\,e_{i-1,1}g_{j,i+1}-m\,g_{j,i+1}e_{i,1}\,\;\;\;\;\;\text{by $(10)$
and $(1)$}\end{eqnarray*} And for the type
\setcounter{c}{2}$(\Roman{c})$:
$$g_{i-1}.b=g_ie_{i-1,1}+m\,e_{i-1,1}-m\,e_{i,1}\qquad\qquad\;\,\,\text{by $(10)$}$$ Thus, we get:
$$g_{i-1}.(b\otimes_{_{H}}1)= \left\lb\begin{array}{l} g_{j,i}e_{i-1,1}\otimes_{_{H}}1+
m\,r^{j-i}\,e_{i-1,1}\ot 1-m\,g_{j,i+1}e_{i,1}\ot1\\
g_ie_{i-1,1}\ot1+m\,e_{i-1,1}\ot1-m\,e_{i,1}\ot1
\end{array}\right.$$
where the first line refers to type \setcounter{c}{1}$(\Roman{c})$
and the second line to type \setcounter{c}{2}$(\Roman{c})$. For
future references, we name these two equalities
\setcounter{c}{1}$A(\Roman{c})$ and \setcounter{c}{2}$A(\Roman{c})$
respectively.
\subsubsection{Action by $g_i$ (Case $B$)}
We have for types \setcounter{c}{1}(\Roman{c}) and
\setcounter{c}{2}(\Roman{c}) respectively:
\begin{eqnarray*}g_i.b&=&g_{j,i+2}e_{i+1,1}\qquad\;\;\;\;\;\;\;\;\;\;\;\;\;\;\;\text{by
$(9)$}\\
g_i.b&=&l^{-1}\,e_{i,1}\qquad\qquad\;\;\;\;\;\;\;\;\;\;\;\;\,\;\;\text{by
$(4)$}\end{eqnarray*} Thus, we get:
$$g_i.(b\otimes_{_{H}}1)=\left\lb\begin{array}{l} g_{j,i+2}e_{i+1,1}\otimes_{_{H}}1\;\;\;\;\setcounter{c}{1}B(\Roman{c})\\
l^{-1}\,e_{i,1}\ot1 \qquad\,\,\;\; \setcounter{c}{2}B(\Roman{c})
\end{array}\right.$$
Notice that if $j=i+1$, expression \setcounter{c}{1}$B(\Roman{c})$
reduces to $e_{i+1,1}\ot1$. \subsubsection{Action by $g_{j}$ (Case
$C$)} Let's first deal with Type \setcounter{c}{1}(\Roman{c}). We
have by $(7)$:
\begin{equation*}g_j.b=g_{j-1,i+1}e_{i,1}-m\,g_{j,i+1}e_{i,1}
+ml^{-1}\,e_jg_{j-1,i+1}e_{i,1}\end{equation*} We will rearrange the
last term of the sum and to do so, we will need more mixed braid
relations, as in the following lemma.
\newtheorem{Lemma}{Lemma}
\begin{Lemma}
\begin{eqnarray}
e_je_ig_j&=&e_jg_i^{-1}\;\;\;\text{when $i\sim j$}\\
e_ie_je_i&=&e_i\;\;\;\;\;\;\;\;\;\text{when $i\sim j$}
\end{eqnarray}
\end{Lemma}
\noin\noin\textbf{Proof.} These are equalities $(8)$ and $(10)$ of
Proposition $2.3$ in \cite{CGW}. $\;\;\;\square$ \\\\
Using the relations of Lemma $1$, we now give a new expression for
$e_jg_{j-1,i+1}e_{i,1}$.
\begin{Lemma}
\begin{equation}e_jg_{j-1,i+1}e_{i,1}=e_{j,1}g_j^{-1}\dots
g_{i+2}^{-1}
\end{equation}
\end{Lemma}
\noin\noin\textbf{Proof.} Using Lemma $1$, we replace $e_j$ with
$e_je_{j-1}e_j$, then replace $e_{j-1}e_jg_{j-1}$ with
$e_{j-1}g_j^{-1}$ to get:
$$e_jg_{j-1,i+1}e_{i,1}=e_je_{j-1}g_{j-2,i+1}e_{i,1}g_j^{-1}$$
Proceeding inductively, we obtain $(13)$. $\;\;\;\square$ \\\\
When $b$ is of the second type, we simply have:
$$g_j.b=\begin{cases}
e_{i,1}g_j&\text{if $j>i+1$}\\g_{i+1}e_{i,1}&\text{if
$j=1+1$}\end{cases}$$ \noin In the following expressions, the first
line is for type \setcounter{c}{1} $(\Roman{c})$ and the next two
are for type \setcounter{c}{2} $(\Roman{c})$.
$$g_j.(b\otimes_{_{H}}1)=\left\lb\begin{array}{l} g_{j-1,i+1}e_{i,1}\ot
1-m\,g_{j,i+1}e_{i,1}\ot 1 +\frac{m}{lr^{j-i-1}}\,e_{j,1}\ot1\\
r\,e_{i,1}\ot 1\;\;\;\;\;\,\,\!\;\;\;\text{when $j>i+1$}\\
g_{i+1}e_{i,1}\ot 1\;\;\!\;\;\;\text{when $j=i+1$}
\end{array}\right.$$
We will refer to these equations as \setcounter{c}{1}$(C\Roman{c})$,
\setcounter{c}{2}$(C\Roman{c})$ and $(C\Roman{c}^{\;'})$
respectively. \subsubsection{Action by $g_{j+1}$ (Case $D$)} Since
\begin{eqnarray*}
g_{j+1}.b&=&g_{j+1}g_{j,i+1}e_{i,1}\qquad\text{when $b$
is of type \setcounter{c}{1}$(\Roman{c})$}\\
g_{j+1}.b&=&e_{i,1}g_{j+1}\qquad\qquad\,\,\text{when $b$ is of type
\setcounter{c}{2}$(\Roman{c})$,}
\end{eqnarray*}
we simply get:
$$g_{j+1}.b=\left\lb\begin{array}{l}
g_{j+1,i+1}e_{i,1}\ot
1\qquad\setcounter{c}{1}(D\Roman{c})\\
r\,e_{i,1}\ot 1\qquad\qquad\;\,\;\setcounter{c}{2}(D\Roman{c})
\end{array}\right.$$
\subsubsection{Action by $g_k$ where $k\not\in\lb i-1,i,j,j+1\rb$
(Case $E$)} \indent $\;\;\;\;\bullet\;\;$ Suppose first $k<i-1$ and
$b$ of type \setcounter{c}{1}(\Roman{c}). We compute:
$$g_k.b=g_{j,i+1}e_{i,k+2}g_ke_{k+1}e_{k,1}$$
Expanding $g_ke_{k+1}e_k$ with $(10)$ yields:
\begin{equation*}
g_k.b=g_{j,i+1}e_{i,k+2}g_{k+1}e_{k,1}+m\,g_{j,i+1}e_{i,k+2}e_{k,1}-m\,g_{j,i+1}e_{i,1}
\end{equation*} Since $e_{k+2}e_k=0$ in $B_1$, this expression
simplifies as follows:
$$g_k.b=g_{j,i+1}e_{i,k+2}(g_{k+1}-m\,e_{k+1})e_{k,1}$$
Replacing $g_{k+1}-m\,e_{k+1}=g_{k+1}^{-1}-m$ with $(8)$ and
simplifying with $e_{k+2}e_k=0$, we then obtain:
$$g_k.b=g_{j,i+1}e_{i,k+2}g_{k+1}^{-1}e_{k,1}$$
Applying equality $(11)$ to $e_{k+2}g_{k+1}^{-1}$ now yields the new
expression for $g_k.b$:
$$g_k.b=g_{j,i+1}e_{i,k+1}g_{k+2}e_{k,1},$$
which is also after commutation of $g_{k+2}$:
$$g_k.b=g_{j,i+1}e_{i,1}g_{k+2}$$
Thus, the action of $g_k$ on the tensor $g_{j,i+1}e_{i,1}\ot 1$ is simply a multiplication by $r$. After inspection, the computation for type \setcounter{c}{2}(\Roman{c}) is identical and the action of $g_k$ on the tensor $e_{i,1}\ot 1$ is also a multiplication by $r$. \\\\
\indent $\bullet\;\;$ Suppose now that $k>j+1$. Visibly, $g_k$
commutes with $g_{j,i+1}e_{i,1}$ and with $e_{i,1}$,
so that in both cases, the action by $g_k$ on the tensor $b\ot 1$ is simply a multiplication by $r$.\\\\
\indent $\bullet\;\;$ Finally, suppose $k$ belongs to $\lb
i+1,\dots,j-1\rb$ where $j\geq i+2$. We look at the action of $g_k$
on $g_{j,i+1}e_{i,1}$. We move $g_k$ to the right, then use the
Braid relation $g_kg_{k+1}g_k=g_{k+1}g_kg_{k+1}$, then move
$g_{k+1}$ this time to the end of the expression. After doing these
moves, we get:
$$g_kg_{j,i+1}e_{i,1}=g_{j,i+1}e_{i,1}g_{k+1}$$
It follows in particular that $$g_k.\,b\ot 1=r\,b\ot 1,$$ as in the
previous cases. It remains to look at the action of $g_k$ on
$e_{i,1}$. We have:
$$g_k.b=\begin{cases} g_{i+1}e_{i,1}&\text{if $k=i+1$}\\
e_{i,1} g_{k}&\text{otherwise}\end{cases}$$ We summarize Case $E$ in
the following two equalities:
$$\forall k\not\in\lb i-1,i,j,j+1\rb,\, g_k.\,b\,\ot 1=\begin{cases} g_{i+1}b\ot 1&\text{if $k=i+1$ and $b$ of type \setcounter{c}{2}$(\Roman{c})$}\\
r\,b\ot1 &\text{in all the other cases}\end{cases}$$ We note that
the top equality is \setcounter{c}{2}$(C\Roman{c}^{\,'})$. Let's
name the bottom equality $(E)$.

With Cases $A$, $B$, $C$, $D$, $E$, the action of the $g_i$'s on the
vector space $B_1\ot F$ is entirely described. The object of the
next part is to give an expression of the representation in terms of
roots.
\subsection{Expression of the representation in the Lawrence-Krammer
space}
\subsubsection{The Lawrence-Krammer representation}
Following our discussion at the end of $\S\,1.3.2$ and in $\S\,2.1$,
there is a bijection:
$$u:\left.\begin{array}{ccc}
\phi^{+}& \lra & \B\\
\be & \longmapsto & b\end{array},\right.$$ where $b$ is the algebra
element associated with the positive root $\be$, as in $\S\,1.3.2$.
It follows that there is a natural isomorphism $\varphi$ of vector
spaces over $F$, defined on the basis vectors by:
$$\varphi:\left.\begin{array}{ccc}
\V^{(n)}& \widetilde{\lra} & B_1\ot F\\
x_{\be} & \longmapsto &  u(\be)\ot 1\end{array}\right.$$ We now get
a representation of the BMW algebra inside the Lawrence-Krammer
space as follows.
\begin{theo}
The map on the generators
$$\n^{(n)}:\left.\begin{array}{ccc}
B(A_{n-1})& \lra & End_{F}(\V^{(n)})\\
g_i & \longmapsto & \n_i\end{array},\right.$$ where each $\n_i$ is
defined on the basis vectors of $\V^{(n)}$ by
$$\n_i(x_{\be})= \varphi^{-1}(g_i.(u(\be)\ot 1))$$
defines a representation of degree $\chl$ of the BMW algebra
$B(A_{n-1})$ in the L-K space $\V^{(n)}$. Once irreducibility over
$\Q(l,r)$ has been established, as a representation of the Braid
group on $n$ strands, it is equivalent to the Lawrence-Krammer
representation.
\end{theo}
\noin\textbf{Proof.} By definition, $\n^{(n)}$ is a representation
of $B$ in $\V^{(n)}$. We notice that $\n^{(n)}$ factors through the
quotient $B/I_2$, where $I_2$ is the two sided ideal of $B$
generated by all the products $e_ie_j$ with $|i-j|>2$. Indeed, in
$B_1$, the algebra element $e_ie_j$ is zero and so in $B_1\ot F$,
the vector $e_ie_jb\ot 1$ is zero. Thus, we have:
$$\n^{(n)}(e_ie_j)=0\;\;\text{when $|i-j|>2$}$$
Then by \cite{CGW}, as a representation of the Braid group on $n$
strands, $\n^{(n)}$ must be equivalent to the Lawrence-Krammer
representation of the Artin group of type $A_{n-1}$ based on the two
parameters $t$ and $r$, as defined in \cite{CW}.
The $r$ of this paper is the $\unsurr$ of \cite{CW}; the parameter
$t$ of $\cite{CW}$ is related to the parameters $l$ and $r$ of this
paper by $lt=r^3$. $\;\;\;\square$

\subsubsection{An explicit form of the representation in terms of
roots} Given a positive root $$\be=\al_i+\dots\al_j\;\;\text{with
$i<j$},$$ we read on the expressions \setcounter{c}{1}
$(A\Roman{c})$, $(B\Roman{c})$, $(C\Roman{c})$, $(D\Roman{c})$ and
$(E)$ of $\S2.2$ an expression of $\n_k(x_{\be})$ for $k\in\lb
i-1,i,j,j+1\rb$ and for $k\not\in\lb i-1,i,j,j+1\rb$. We define the
height $ht(\be)$ of a positive root $\be$ as the sum of its
coefficients with respect to the simple roots
$\al_1,\dots,\al_{n-1}$. We have:
\begin{eqnarray*}
\n_{i-1}(x_{\be})&=&x_{\be+\al_{i-1}}+m\,r^{ht(\be)-1}\,x_{\al_{i-1}}-m\,x_{\be}\;\;\;\;\!\,\;\;\,\text{by
$(A\Roman{c})$}\\
\n_i(x_{\be})&=&x_{\be-\ali}\qquad\qquad\qquad\qquad\qquad\qquad\;\;\;\;\,\!\;\;\;\,\text{by
$(B\Roman{c})$}\\
\n_j(x_{\be})&=&x_{\be-\al_j}+\frac{m}{l\,r^{ht(\be)-2}}\,x_{\al_j}-m\,x_{\be}\;\;\;\,\;\;\;\;\;\;\;\;\;\;\,\!\text{by $(C\Roman{c})$}\\
\n_{j+1}(x_{\be})&=&x_{\be+\al_{j+1}}\qquad\qquad\qquad\qquad\qquad\;\;\;\;\;\;\;\;\;\,\;\;\,\!\text{by
$(D\Roman{c})$}\\
\n_k(x_{\be})&=&r\,x_{\be}\qquad\forall k\not\in\lb
i-1,i,j,j+1\rb\;\;\;\;\;\;\,\;\,\,\!\,\;\;\,\text{by $(E)$}
\end{eqnarray*}
Similarly for type \setcounter{c}{2}$(\Roman{c})$, if $\be=\al_i$ is
a simple root, we have:
\begin{eqnarray*}
\n_{i-1}(\xb)&=& x_{\be+\al_{i-1}}+m\,x_{\al_{i-1}}-m\,\xali\;\;\;\;\!\,\;\;\,\;\;\;\;\;\;\;\;\;\;\;\;\;\text{by $(A\Roman{c})$}\\
\n_i(x_{\be})&=& l^{-1}\,\xali\qquad\qquad\qquad\qquad\qquad\qquad\;\;\;\;\,\!\;\;\;\,\text{by $(B\Roman{c})$}\\
\n_{i+1}(x_{\be})&=& x_{\be+\al_{i+1}}\qquad\qquad\qquad\qquad\qquad\;\;\;\;\;\;\;\;\;\,\;\;\,\!\;\;\text{by $(C\Roman{c}^{\,'})$}\\
\n_k(x_{\be})&=&r\,\xb\qquad\forall k\not\in\lb i-1,i,i+1\rb
\end{eqnarray*}
The last equation is obtained with $(C\Roman{c})$ when $k>i+1$ and
with $(E)$ when $k<i-1$.\\
For each node $i$, we summarize the action of $\n_i$ on $\xb$ as
follows.
$$\n_i(\xb)=\begin{cases}
r\,\xb&\text{if $(\be|\al_i)=0\;\;\;(a)$}\\
\unsur{l}\,\xb&\text{if $(\be|\al_i)=1\;\;\;(b)$}\\
x_{\be+\al_i}&\text{if $(\be|\al_i)=-\unsur{2}\;\;$ and (c)}\\
x_{\be+\al_i}+m\,r^{ht(\be)-1}\,\xali-m\,\xb&\text{if
$(\be|\al_i)=-\unsur{2}\;\;$ and $(c^{'})$}\\
x_{\be-\al_i}+\frac{m}{l\,r^{ht(\be)-2}}\,\xali-m\,\xb&\text{if
$(\be|\al_i)=\unsur{2}\;\;\;\;\;$ and $(d)$}\\
x_{\be-\al_i}&\text{if $(\be|\al_i)=\unsur{2}\;\;\;\;\;$ and
$(d^{'})$}
\end{cases}$$

\noin where $(c), (c^{'}), (d),(d^{'})$ are the following
conditions: $$\begin{array}{cccccc}
(c)& \be&=&\al_t+\dots+\al_{i-1}&\text{with}&t\leq i-1\\
(c^{'})& \be&=&\al_{i+1}+\dots +\al_s&\text{with}&s\geq i+1\\
(d)&\be&=&\al_t+\dots+\ali&\text{with}&t\leq i-1\\
(d^{'})& \be&=&\ali+\dots+\al_s&\text{with}&s\geq i+1
\end{array}$$
We note that: \\$(a)$ is equivalent to $Supp(\be)\cap\lb
i-1,i,i+1\rb=\emptyset$ or
$\lb i-1,i,i+1\rb\subseteq Supp(\be)$. \\
$(b)$ is equivalent to $\be=\ali$. \\
$(c)$ and $(c^{'})$ are the two ways the inner product
$(\be|\ali)$ can be $\unsur{2}$.\\
$(d)$ and $(d^{'})$ are the two ways the inner product $(\be|\ali)$
can be $-\unsur{2}$.\\\\ We deduce from these equalities an
expression for $\n^{(n)}(e_i)$:
$$\n^{(n)}(e_i)(\xb)=\begin{cases}
0&\text{if $(\be|\al_i)=0$}\\
\bigg(1-\frac{l-\unsur{l}}{\unsurr-r}\bigg)\;\xali&\text{if $(\be|\al_i)=1$}\\
\unsur{r^{ht(\be)-1}}\;\;\;\;\,\;\;\;x_{\al_i}&\text{if $(\be|\al_i)=-\unsur{2}\;\;$ and (c)}\\
r^{ht(\be)-1}\;\;\;\;\;\;\,\xali &\text{if
$(\be|\al_i)=-\unsur{2}\;\;$ and $(c^{'})$}\\
\unsur{l\,r^{ht(\be)-2}}\;\;\;\,\!\;\;\;\xali&\text{if
$(\be|\al_i)=\unsur{2}\;\;\;\;\;$ and $(d)$}\\
l\,r^{ht(\be)-2}\;\;\;\;\;\!\xali&\text{if
$(\be|\al_i)=\unsur{2}\;\;\;\;\;$ and $(d^{'})$}
\end{cases}$$

\noin Notice $\n^{(n)}(e_i)(\xb)$ is always a multiple of
$\xali$. This is easily pictured on the tangles.\\
The next part establishes Theorem $1$, following the discussion of
$\S 1.2$.

\section{Reducibility of the representation}
\subsection{Action on a proper invariant subspace of the L-K space}
We show the following result:
\newtheorem{Proposition}{Proposition}
\begin{Proposition}\hfill\\
For any proper invariant subspace $\U$ of $\V^{(n)}$, we have
$\n^{(n)}(e_i)(\U)=0$ for all $i$.
\end{Proposition}
\noin\textbf{Proof.} If $\U$ is trivial, there is nothing to prove.
Otherwise, let $u$ be a nonzero vector of $\U$ such that
$\n^{(n)}(e_i)(u)\neq 0$. Since $\n^{(n)}(e_i)(u)$ is a multiple of
$\xali$, we see that $\xali$ is in $\U$. 
From there, we have:
$$\n_{i-1}(\xali)=x_{\al_i+\al_{i-1}}+m\,x_{\al_{i-1}}\;\;\;modulo\;\;\;F\xali$$
Hence $x_{\al_i+\al_{i-1}}+m\,x_{\al_{i-1}}$ is in $\U$. Another
application of $\n_{i-1}$ now yields:
$$\n_{i-1}(x_{\al_i+\al_{i-1}})+m\,x_{\al_{i-1}}=\xali+\frac{m}{l}\,x_{\al_{i-1}},$$
from which we derive that $x_{\al_{i-1}}$ is in $\U$. By induction,
we see that all the $x_{\al_t}$'s for $t\leq i$ are in $\U$. In
particular, $x_{\al_1}$ is in $\U$. But since $e_1\ot 1$ spans
$B_1\ot F$, $x_{\al_1}$ spans $\V^{(n)}$. Then $\U$ is the whole L-K
space $\V^{(n)}$, in contradiction with our assumption that $\U$ is
proper.  $\;\;\;\square$

\newtheorem{Corollary}{Corollary}
\begin{Corollary}
Let $\W$ be a proper irreducible invariant subspace of $\V^{(n)}$.
Then, $\W$ is an irreducible $\ih(n)$-module.
\end{Corollary}

\noin\textbf{Proof.} By Proposition $1$ and $(3)$, we have
$$\big[g_i^2+m\,g_i-1\big].\W=0\;\;\;\text{for all $i$.}$$
Hence $\W$ is an $\ih(n)$-module. Since the $e_i$'s are polynomials
in the $g_i$'s, $\W$ is an irreducible $\ih(n)$-module.
$\;\;\;\square$\\

The next part investigates the existence of a one-dimensional
invariant subspace of $\V^{(n)}$. We define for two nodes $i$ and
$j$ with $i<j$
$$w_{ij}=x_{\ali+\dots+\al_{j-1}}$$
We will sometimes write $w_{i,j}$ instead of $w_{ij}$. Below is how
$w_{ij}$ is represented in the tangle algebra:\\\vspace{-0.28in}
\begin{center}
\epsfig{file=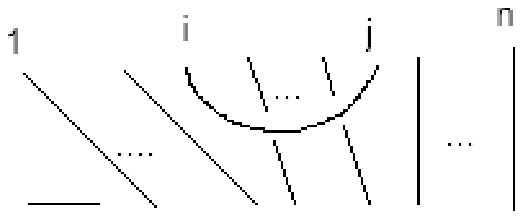, height=4cm}
\end{center} \vspace{-0.5in}

It has two horizontal strands: one that joins nodes $i$ and $j$ at
the top, and one that joins nodes $1$ and $2$ at the bottom and
$(n-2)$ vertical strands that don't cross within each other. The top
horizontal strand over-crosses the vertical strands that it
intersects.
\subsection{\textbf{The case $l=\unsur{r^{2n-3}}$}}

We will prove the theorem:
\begin{theo} Let $n$ be an integer with $n\geq 3$ and assume $(r^2)^2\neq 1$.\\
\indent Suppose $n\geq 4$. There exists a one-dimensional invariant
subspace of $\V^{(n)}$ if and only if $l=\unsur{r^{2n-3}}$. If so,
it is spanned by $\sum_{1\leq s<t\leq
n}r^{s+t}\,w_{st}$\\
\indent\textit{(Case $n=3$)} There exists a one-dimensional
invariant subspace of $\V^{(3)}$ if and only if $l=\unsur{r^3}$ or
$l=-r^3$.\vspace{0.05in} \\Moreover, if $r^6\neq -1$, it is unique
and \vspace{-0.099in}
$$\begin{array}{cccccc} \text{when}&l&=&\unsur{r^3},&\text{it is
spanned
by}&w_{12}+r\,w_{13}+r^2\,w_{23}\\
\text{when}&l&=&-r^3,&\text{it is spanned
by}&w_{12}-\unsurr\,w_{13}+\unsur{r^2}\,w_{23}
\end{array}\vspace{-0.11in}$$
\noin If $r^6=-1$, there are exactly two one-dimensional invariant
subspaces of $\V^{(3)}$ and they are respectively spanned by the
vectors above.
\end{theo}
\noin\textbf{Proof.} Let $\U$ be a one-dimensional invariant
subspace of $\V^{(n)}$ and $u$ a spanning vector of $\U$. For each
$i$, let $\gamma_i$ be the scalar such that $\n_i(u)=\gamma_i\,u$.
Since $(\n_i^2+m\,\n_i-id_{\V^{(n)}})(u)=0$ by Proposition $1$, it
follows that $\gamma_i^2+m\,\gamma_i-1=0$. Hence $\gamma_i\in\lb
r,-\unsurr\rb$. Further, since $(r^2)^2\neq 1$, the Braid relation
$\n_i\n_j\n_i=\n_j\n_i\n_j$ when $i\sim j$ forces that $\gamma_i$
takes the same value as $\gamma_j$. Let's denote by $\gamma$ the
common value of the $\gamma_i$'s. So, for each $i$, we have
$\n_i(u)=\gamma\,u$, where $\gamma\in\lb r,-\unsurr\rb$. \\
A general form for $u$ is:
$$u=\sum_{1\leq i<j\leq n}\mu_{ij}\,w_{ij},\;\;\;\text{where $\;\;\mu_{ij}\in F$}$$
We look for relations between these coefficients.
\begin{Lemma}
Let $i$ be some node. Suppose $v=\sum_{1\leq k<f\leq
n}\mu_{kf}\,w_{kf}$ is a vector of $\V^{(n)}$ with
$\n_i(v)=\gamma\,v$ where $\gamma\in\lb r,-\unsurr\rb$. Then the
following equalities hold for the coefficients of $v$:
\begin{eqnarray}
\forall s\geq i+2,\,\mu_{i+1,s}&=&\gamma\,\mu_{i,s}\\
\forall t\leq i-1,\,\mu_{t,i+1}&=&\gamma\,\mu_{t,i}
\end{eqnarray}
When $i=1$, only $(14)$ holds and when $i=n-1$, only $(15)$ holds.
\end{Lemma}
\noin\textbf{Proof.} To show $(14)$, we look at the coefficient of
$w_{i+1,s}$ in $\n_i(v)=\gamma\,v$, where $s\geq i+2$. We get:
$\mu_{i,s}-m\,\mu_{i+1,s}=\gamma\,\mu_{i+1,s}$. Since
$\gamma+m=\unsur{\gamma}$, this equality is equivalent to
$\mu_{i+1,s}=\gamma\mu_{i,s}$. Similarly, by equating the
coefficients of $w_{t,i+1}$ in $\n_i(v)=\gamma\,v$, we obtain
$(15)$.
$\;\;\;\square$\\\\
Applying these equalities to the coefficients of $u$, we see that
all the coefficients of $u$ must be nonzero. In particular, when
$n\geq 4$, the coefficient $\mu_{34}$ of $u$ is nonzero. Because an
action of $g_1$ on $w_{34}$ is a multiplication by $r$ and an action
on $g_1$ on the other terms $w_{ij}$ does not create any term in
$w_{34}$, this forces $\gamma=r$. Thus, without loss of generality,
we have:
$$u=\sum_{1\leq i<j\leq n} r^{i+j}\,w_{ij}$$
From there, let's look at the action of $g_1$ on $u$ and the
resulting coefficient in $w_{12}$. The action of $g_1$ on $w_{12}$
is a multiplication by $l^{-1}$ and an action of $g_1$ on the
$w_{2,j}$'s for $3\leq j\leq n$ creates new terms in $w_{12}$ with
respective coefficients $m\,r^{j-3}$. Thus, we get the equation:
\begin{equation*}
\frac{r^3}{l}+\sum_{j=3}^n(r^2)^j=r^4,
\end{equation*}
from which we derive that $l=\unsur{r^{2n-3}}$.\\
Conversely, if $l=\unsur{r^{2n-3}}$, we define $u$ as $\sum_{1\leq
i<j\leq n}r^{i+j}\,w_{ij}$ and check that $\n_i(u)=r\,u$ for each
$i$. First we show that the coefficient $r^{2i+1}$ of $w_{i,i+1}$ is
multiplied by $r$ when acting by $g_i$. There are three
contributions. One comes from the terms $w_{t,i+1}$, with
coefficient $r^{t+i+1}\,\frac{m\,r^{2n-3}}{r^{i-t-1}}$. Another one
from the terms $w_{i+1,s}$ with coefficient
$r^{s+i+1}\,m\,r^{s-i-2}$, and a third one simply from the term
$w_{i,i+1}$, with coefficient $r^{2i+1}r^{2n-3}$. Now, we have:
$$mr^{2n-1}\sum_{t=1}^{i-1}r^{2t}+\frac{m}{r}\sum_{s=i+2}^n
r^{2s}+r^{2n+2i-2}=r\,r^{2i+1}$$ Next, given a positive root $\be$,
if none of the nodes $i-1,\,i,\,i+1$ is in the support of $\be$ or
if all three nodes $i-1,\,i,\,i+1$ are in the support of $\be$, then
it comes:
$$\n_i(x_{\be})=r\,x_{\be}$$
Thus, we only need to study the action of $\n_i$ on $w_{k,i}$,
$w_{k,i+1}$, with $k\leq i-1$ on one hand and $w_{i,l}$,
$w_{i+1,l}$, with $l\geq i+2$ on the other hand. \\We have:
\begin{eqnarray*}
r^{k+i}\;\n_i(w_{k,i})&=&r^{k+i}\;w_{k,i+1}\\
r^{k+i+1}\;\n_i(w_{k,i+1})&=&r^{k+i+1}\;w_{k,i}-m\,r^{k+i+1}\;w_{k,i+1}\;
modulo\, F\,\xali
\end{eqnarray*}
So we get:
$$\n_i\Big(r^{k+i}\,w_{k,i}+r^{k+i+1}\,w_{k,i+1}\Big)=r^{k+i+1}\,w_{k,i}+r^{k+i+2}\,w_{k,i+1}\;modulo\,F\,\xali$$
Similarly, we have:
\begin{eqnarray*}
r^{l+i}\;\n_i(w_{i,l})&=&r^{l+i}\;w_{i+1,l}\\
r^{l+i+1}\;\n_i(w_{i+1,l})&=&r^{l+i+1}\;w_{i,l}-m\,r^{l+i+1}\;w_{i+1,l}\;modulo\,
F\,\xali,
\end{eqnarray*}

\noindent so that:
$$\n_i\Big(r^{l+i}\,w_{i,l}+r^{l+i+1}\,w_{i+1,l}\Big)=r^{l+i+1}\,w_{i,l}+r^{l+i+2}\,w_{i+1,l}\;modulo\,F\,\xali$$
This ends the proof of the Theorem when $n\geq 4$.\\ Suppose now
$n=3$. So, $$u=\gamma^3\,w_{12}+\gamma^4\,w_{13}+\gamma^5\,w_{23}$$
Let's compute $\n_1(u)$ and $\n_2(u)$:
$$\begin{array}{ccccccccc}
\n_1(u) & = & (\frac{\gamma^3}{l}+\,m\,\gamma^5) & w_{12} & + &
\gamma^5\,w_{13} & + & \gamma^6 & w_{23}\\
\n_2(u) & = & \gamma^4 & w_{12} & + & \gamma^5\,w_{13} & + &
(\frac{\gamma^5}{l}+\frac{m\,\gamma^4}{l}) & w_{23}
\end{array}$$
Since $\n_1(u)=\gamma\,u$, we must have:
$$\begin{array}{ccccccc} \frac{\gamma^3}{l}+m\,\gamma^5=\gamma^4,&
\textit{i.e.}& \unsur{l}=\gamma(1-m\,\gamma),&\textit{i.e.}&
l=\unsur{\gamma^3}, &\text{as}& 1-m\,\gamma=\gamma^2\end{array}$$
\noindent Thus, if there exists a one dimensional invariant subspace
of $\V^{(3)}$ then $l$ must take the values $\unsur{r^3}$ or
$-r^3$.\\ Conversely, let's consider the two vectors:
$$\begin{array}{cccccccc}
&u_r &=& w_{12} &+& r\,w_{13} &+& r^2\,w_{23}\\
&u_{-\unsurr} &=& w_{12} &-& \unsurr\,w_{13}&+&\unsur{r^2}\,w_{23}
\end{array}$$
We read on the equations giving the expressions for $\n_1(u)$ and
$\n_2(u)$ that
$$\begin{array}{cccccccccc} \text{If} &
l&=&\unsur{r^3}&\text{then} & \n_1(u_r)&=&\n_2(u_r)&=&r\;u_r\\
\text{If} & l &=& -r^3 & \text{then} &
\n_1(u_{-\unsurr})&=&\n_2(u_{-\unsurr})&=&-\unsurr\;u_{-\unsurr}
\qquad\square\end{array}$$ The next section investigates the
existence of an irreducible $(n-1)$-dimensional invariant subspace
of $\V^{(n)}$.
\subsection{\textbf{The case
$l\in\lb\unsur{r^{n-3}},-\unsur{r^{n-3}}\rb$}} In Theorem $4$, the
case $n=3$ was special. Likewise, in the following Theorem $5$, the
case $n=4$ needs to be formulated separately.
\begin{theo}\hfill\\
Let $n$ be a positive integer with $n\geq 3$ and $n\neq 4$. Assume
$\ih(n)$ is semisimple. Then, there exists an irreducible
$(n-1)$-dimensional invariant subspace of $\V^{(n)}$ if and only if
$l=\unsur{r^{n-3}}$ or
$l=-\unsur{r^{n-3}}$. \\\\
If so, it is spanned by the $v_i^{(n)}$'s, $1\leq i\leq n-1$, where
$v_i^{(n)}$ is defined by the formula:
\begin{equation*}\begin{split}v_i^{(n)}=\Big(\unsur{r}\,-\,\unsur{l}\Big)w_{i,i+1}&+\sum_{s=i+2}^{n}r^{s-i-2}(w_{i,s}\,-\,\unsurr\;w_{i+1,s})\\
&+\;\e_l\,\sum_{t=1}^{i-1}r^{n-i-2+t}(w_{t,i}\,-\,\unsurr\;w_{t,i+1})\end{split}\end{equation*}
\begin{center} with \hspace{.1in}
$\begin{cases}\e_{\unsur{r^{n-3}}}\;\;=\,1\\\e_{-\unsur{r^{n-3}}}=-1\end{cases}$
\end{center}
\noin\textit{(Case $n=4$)} Assume $\ih(4)$ is semisimple. Then,
there exists an irreducible $3$-dimensional invariant subspace of
$\V^{(4)}$ if and only if
$l\in\lb\unsurr,-\unsurr,-r^3\rb$.\\\\
If $l\in\lb -\unsurr,\unsurr\rb$, it is spanned by $v_1^{(4)}$, $v_2^{(4)}$, $v_3^{(4)}$\\\\
If $l=-r^3$, it is spanned by the vectors:
$$\left\lb\begin{array}{cccccccccccccc} u_1\eg r\,w_{23}+
w_{13}+(\unsurr+\unsur{r^3})w_{34}-
w_{24}- \unsurr\,w_{14}\\
u_2\eg -r\,w_{12}-r^2\,w_{13}-
\unsurr\,w_{34}-\unsur{r^2}\,w_{24}+(r+\unsurr)\,w_{14}\\
u_3\eg (r+\unsur{r^3})\,w_{12}+\unsurr\,w_{23}- w_{13}+w_{24}-
r\,w_{14}
\end{array}\right.$$
\end{theo}
\noin\textbf{Proof.} Suppose that there exists an irreducible
$(n-1)$-dimensional invariant subspace $\U$ of $\V^{(n)}$.
\newtheorem{Claim}{Claim}
\begin{Claim}
Except in the case when $n=6$, there exists a basis
$(v_1,\dots,v_{n-1})$ of $\U$ such that one of the following two
sets of relations holds:

$$\left.\begin{array}{l}(\tr)\left|\begin{array}{l}\n_t(v_i)=r\;v_i\;\,\!\;\;\;\;\;\;\;\;\;\;\,\,\qquad\qquad\forall\;
t\not\in\{
i-1,\,i,\,i+1\}\\
\n_i(v_i)=-\unsurr\;v_i\qquad\qquad\qquad\!\;\,\!\,\,\forall 1\leq i\leq n-1\\
\n_{i+1}(v_i)=r(v_i+v_{i+1})\;\;\;\;\;\;\;\;\,\;\!\,\,\,\forall 1\leq i\leq n-2\\
\n_{i-1}(v_i)=r\;v_i+\unsurr\;v_{i-1}\;\;\;\;\;\;\;\;\!\,\,\,\forall
2\leq i\leq n-1
\end{array}\right.\\\\
(\bigtriangledown)\left|\begin{array}{l}
\n_t(v_i)=-1/r\;v_i\;\;\;\;\;\;\;\,\!\qquad\qquad\forall\;
t\not\in\{
i-1,\,i,\,i+1\}\\
\n_i(v_i)=r\;v_i\;\;\;\qquad\qquad\;\;\;\;\;\;\;\;\;\;\,\forall 1\leq i\leq n-1\\
\n_{i+1}(v_i)=-1/r(v_i+v_{i+1})\;\;\;\;\;\,\forall 1\leq i\leq n-2\\
\n_{i-1}(v_i)=-1/r\;v_i-r\;v_{i-1}\;\;\;\;\;\!\,\forall 2\leq i\leq
n-1
\end{array}\right.\end{array}\right.$$
\end{Claim}
\noin\textbf{Proof.} We first recall some general fact about the
irreducible representations of the Iwahori-Hecke algebra of the
symmetric group. The following result was established by James for
the irreducible representations of the symmetric group, but applies
here to the Iwahori-Hecke algebra $\ih(n)$ since we work in
characteristic zero and assumed $\ih(n)$ semisimple. By Theorem $6$,
point $(i)$ in \cite{J}, when the characteristic of the field $F$ is
zero and for $n\geq 7$, an irreducible $F\,Sym(n)$-module is either
one of the Specht modules $S^{(n)}$, $S^{(1^n)}$, $S^{(n-1,1)}$,
$S^{(2,1^{n-2})}$ or has dimension greater than $(n-1)$. The
statement is also true when $n=3$ and $n=5$. When $n=4$, the
statement does not hold as $S^{(2,2)}$ has dimension $2$ and when
$n=6$, the statement also fails since $S^{(3,3)}$ and $S^{(2,2,2)}$
both have dimension $5$. In any case, there are exactly two
inequivalent irreducible representations of $F\,Sym(n)$ of degree
$(n-1)$, except in the case $n=6$, when there are exactly four
inequivalent irreducible representations of $F\,Sym(6)$ of degree
$5$. The same statement holds for $\ih(n)$ when the algebra is
semisimple.

Consider now the set of relations $(\tr)$ (resp
$(\bigtriangledown)$). For each $i$, let $M_i$ (resp $N_i$) be the
matrix of the endomorphism $\n_i$ in the basis
$(v_1,\dots,v_{n-1})$. It is a direct verification that the $M_i$'s
(resp $N_i$'s) satisfy the Braid relations and the relation
$M_i^2+m\,M_i=I_{n-1}$ (resp $N_i^2+m\,N_i=I_{n-1}$) for each $i$,
where $I_{n-1}$ is the identity matrix of size $(n-1)$. Hence the
$M_i$'s (resp the $N_i$'s) yield a matrix representation of $\ih(n)$
of degree $(n-1)$.

To show that these two matrix representations are irreducible,
relying on James'statement above, it suffices to check that there is
no one-dimensional invariant subspace of $F^{n-1}$ when $n\neq 4$
and that there is no one-dimensional or irreducible two-dimensional
invariant subspace of $F^3$ when $n=4$. This is the case if
$r^{2n}\neq 1$ when $n\neq 4$ and if $(r^2)^2\neq 1$ and
$(r^2)^4\neq 1$ when $n=4$.

When $n=3$, the two matrix representations are equivalent. When
$n\geq 4$, they are not: visibly, the matrices of one representation
all have the same trace $-\frac{(n-2)}{r}+r$ and the matrices of the
other one all have the same trace $(n-2)r-\unsurr$. These two values
are distinct when $(r^2)^2\neq 1$ and $n\geq 4$. We conclude that
these are the two inequivalent irreducible representations of
$\ih(n)$ when $n\geq 4$
and $n\neq 6$. $\;\;\;\square$\\
In what follows, we assume $n\geq 4$. We will show that it is
impossible to have the second set of relations, except in the case
$n=4$ when it forces $l=-r^3$. Suppose the $v_i$'s satisfy
$(\bigtriangledown)$. The relation $\n_{n-1}(v_1)=-\unsurr\,v_1$
implies that in $v_1$ there are no terms in $w_{s,t}$ for integers
$s,t\in\{1,\,\dots,\,n-2\}$ such that $s<t$. Hence we may write:
\begin{equation*}
v_1=\sum_{j=1}^{n-2}\m_{j,n-1}w_{j,n-1}+\sum_{j=1}^{n-1}\m_{j,n}w_{j,n}
\end{equation*}
Moreover, the relation $\n_3(v_1)=-\unsurr\,v_1$ implies that there
are no terms in $w_{j,k}$ in $v_1$ for $j\geq 5$. Further, the
relations $$\left\lbrace\begin{array}{ccc}
\n_1(w_{2,n-1}) &=& w_{1,n-1}+mr^{n-4}\,w_{12}-m\,w_{2,n-1}\\
\n_1(w_{2,n}) &=& w_{1,n}+mr^{n-3}\,w_{12}-m\,w_{2,n}
\end{array}\right.$$
\noindent imply that: $mr^{n-4}\m_{2,n-1}+mr^{n-3}\m_{2,n}=0,\;$
\emph{i.e} $\;\m_{2,n}=-\unsurr\,\m_{2,n-1},$ \noindent as there is
no term in $w_{12}$ in $v_1$. Furthermore, an application of $(14)$
with $\g=r$ and $i=1$ yields for $s=n-1$ and $s=n$ respectively:
$\mu_{2,n-1}=r\,\mu_{1,n-1}$ and $\mu_{2,n}=r\,\mu_{1,n}$. So, up to
a multiplication by a scalar,
\begin{multline}
v_1=w_{1,n-1}+r\,w_{2,n-1}-\unsurr\,w_{1,n}-w_{2,n}\\
+\m_{3,n-1}\,w_{3,n-1}+\m_{4,n-1}\,w_{4,n-1}+\m_{3,n}\,w_{3,n}+\m_{4,n}\,w_{4,n}
\end{multline}
\begin{equation}
\text{or}\;\;\qquad
v_1=\m_{3,n-1}\,w_{3,n-1}+\m_{4,n-1}\,w_{4,n-1}+\m_{3,n}\,w_{3,n}+\m_{4,n}\,w_{4,n}
\end{equation}
If $n\geq 5$, the relation $\n_3(v_1)=-\unsurr\,v_1$ implies that
$\mu_{2,n}=0$. Then an expression for $v_1$ is given by $(17)$ and
not $(16)$. Assume $n>5$. Then there is no term in $w_{34}$ in
$v_1$. Then it comes $\mu_{4,n}=-\unsurr\,\mu_{4,n-1}$. Moreover, by
$(14)$ applied with $\g=-\unsurr$ and $i=3$, we have
$\mu_{4,n}=-\unsurr\,\mu_{3,n}$ and
$\mu_{4,n-1}=-\unsurr\,\mu_{3,n-1}$. Further, when $n>5$,
$\n_4(w_{3,n})=r\,w_{3,n}$ and an action of $g_4$ on the other terms
of $v_1$ in $(17)$ won't create any term in $w_{3,n}$. Thus, the
relation $\n_4(v_1)=-\unsurr\,v_1$ forces $\mu_{3,n}=0$. Then, by
the relations previously established, all the coefficients of $v_1$
are actually zero, which is impossible. The case $n=5$ also leads to
a contradiction and details appear in \cite{L}, $\S8.3$. \\
Suppose now $n=4$. By $(16)$ and $(17)$,
$v_1=w_{13}+r\,w_{23}-\unsurr\,w_{14}-w_{24}+\mu_{34}\,w_{34}$ or
$v_1=w_{34}$. Suppose $v_1$ is of the second type. Then, by
$\n_3(v_1)=-\unsurr\,v_1$, we must have $l=-r$. Since
$\n_2(v_1)=-\unsurr(v_1+v_2)$, we must have:
\begin{equation*}
v_2=-r\,w_{24}+(r^2-1)\,w_{23}-r^2\,w_{34}
\end{equation*}
Since
$\left\lb\begin{array}{l}\n_2(w_{34})=w_{24}+m\,w_{23}-m\,w_{34}\\\n_2(w_{23})=-\unsurr\,w_{23}\end{array}\right.$
and since $\n_2(v_2)=r\,v_2$, we get:
$$-m\,r^2-\unsurr(r^2-1)=r(r^2-1),\;\text{which reads $m=0$ after simplification.}$$
As $m$ is nonzero, this is a contradiction. Thus, $v_1$ is of the
first type. Then, denoting by $\la_{ij}$ the coefficient of $w_{ij}$
in $v_2$, we get by looking at the coefficient of $w_{12}$ in
$\n_2(v_1)=-\unsurr\,v_1-\unsurr\,v_2$ that $\la_{12}=-r$. Since by
$(15)$ with $\g=r$ and $i=2$, we have $\la_{13}=r\,w_{12}$, it
follows that $\la_{13}=-r^2$. Next, by looking at the coefficient of
$w_{14}$ in the relation $\n_2(v_1)=-\unsurr\,v_1-\unsurr\,v_2$, we
get:
$$-1=-\frac{\la_{14}}{r}+\unsur{r^2}\qquad i.e.\qquad \la_{14}=r+\unsurr$$
Also, by looking at the coefficient of $w_{24}$ in the same
relation, we obtain:
\begin{equation}\mu_{34}=\unsurr-\frac{\la_{24}}{r}\end{equation}
Next, we use the relation $\n_1(v_2)=-\unsurr\,v_2-r\,v_1$. First we
look at the coefficient of $w_{13}$ to get $\la_{23}=0$ and by
looking at the coefficient of $w_{14}$, we get
$\la_{24}=-\unsur{r^2}$. By $\n_2(v_2)=r\,v_2$ and $(14)$, it comes
$\la_{34}=r\,\la_{24}=-\unsurr$. Plugging the value of $\la_{24}$ in
$(18)$ now yields $\mu_{34}=\unsurr+\unsur{r^3}$. Finally, by
looking at the coefficient of $w_{12}$ in
$\n_1(v_2)=-\unsurr\,v_2-r\,v_1$, we get
$-\frac{r}{l}-\frac{m}{r}=1$, from which we derive $l=-r^3$. Also,
gathering all the results above, we see that the vectors $v_1$ and
$v_2$ are exactly the vectors $u_1$ and $u_2$ of the Theorem.
Similar computations would also lead to $v_3=u_3$ (see \cite{L}).
Conversely, we can show that if $l=-r^3$, then the vectors $u_1$,
$u_2$ and $u_3$ form a free family of vectors that satisfy the
relations $(\bigtriangledown)$. This shows that their linear span
over $F$ is an irreducible $3$-dimensional invariant subspace of
$\V^{(4)}$. For details, see \cite{L}, $\S8.3$.\\
Suppose now that the $v_i$'s satisfy $(\tr)$. The relation
$\n_i(v_i)=-\unsurr\,v_i$ implies that in $v_i$ there are no terms
in $w_{ts}$ for $s\leq i-1$ or $t\geq i+2$ or $t\leq i-1$ and $s\geq
i+2$. Thus, a general form for $v_i$ must be:
\begin{equation}
\begin{split}
v_i=\m_{i,i+1}\,w_{i,i+1} & +\sum_{s=i+2}^n\m_{i,s}\,w_{i,s}+\sum_{s=i+2}^n\m_{i+1,s}\,w_{i+1,s}\\
&\qquad\qquad+\sum_{t=1}^{i-1}\m_{t,i}\,w_{t,i}+\sum_{t=1}^{i-1}\m_{t,i+1}\,w_{t,i+1}
\end{split}
\end{equation}
Since $\n_i(v_i)=-\unsurr\,v_i$, both equalities $(14)$ and $(15)$
hold with $\gamma=-\unsurr$. Further, since $\n_q(v_i)=r\,v_i$ for
$q\not\in\lb i-1,i,i+1\rb$, applying $(14)$ and $(15)$ with $i=q$
and $\gamma=r$ yields:
\begin{eqnarray}
\forall j\geq q+2,\,\mu_{q+1,j}&=&r,\mu_{q,j}\\
\forall k\leq q-1,\,\mu_{k,q+1}&=&r\,\mu_{k,q}
\end{eqnarray}
Apply $(20)$ with $q\leq i-2$ and $j\in\lb i,i+1\rb$ to get:
$$\forall q\leq
i-2,\,\mu_{q+1,i}=r\,\mu_{q,i}\;\;\&\;\;\mu_{q+1,i+1}=r\,\mu_{q,i+1}$$
Apply $(21)$ with $q\geq i+2$ and $k\in\lb i,i+1\rb$ to get:
$$\forall q\geq
i+2,\,\mu_{i,q+1}=r\,\mu_{i,q}\;\;\&\;\;\mu_{i+1,q+1}=r\,\mu_{i+1,q}$$
Expression $(19)$ may now be rewritten:
$$v_i=\zeta^{(i)}\,w_{i,i+1}+\delta^{(i)}\,\sum_{s=i+2}^n
r^{s-i-2}(w_{i,s}-\unsurr\,w_{i+1,s})+\la^{(i)}\,\sum_{t=1}^{i-1}r^{t-1}\,(w_{t,i}-\unsurr\,w_{t,i+1}),$$
where $\zeta^{(i)}$, $\delta^{(i)}$ and $\la^{(i)}$ are three
coefficients to determine. First, we show that all the
$\delta^{(i)}$ with $i\in\lb 1,\dots, n-2\rb$ may be set to the
value one. Notice that if $v_1,\dots,v_{n-1}$ satisfy $(\tr)$, then
$\delta\,v_1,\dots,\delta\,v_{n-1}$ also satisfy $(\tr)$, where
$\delta$ is any nonzero scalar. Then, without loss of generality, we
set $\delta^{(1)}=1$. Suppose $\delta^{(i)}=1$ for some node $i$
with $1\leq i\leq n-2$. We will show that $\delta^{(i+1)}=1$. Notice
that $\delta^{(i+1)}$ is the coefficient of $w_{i+1,i+3}$ in
$v_{i+1}$. Since an action of $g_{i+1}$ on $v_i$ never creates a
term in $w_{i+1,i+3}$, by looking at the coefficient of
$w_{i+1,i+3}$ in $\n_{i+1}(v_i)=r\,v_i+r\,v_{i+1}$, we get
$0=-r\,\delta^{(i)}+r\,\delta^{(i+1)}$. After replacing
$\delta^{(i)}$ by $1$, this yields $\delta^{(i+1)}=1$. Thus, all the
$\delta^{(i)}$ may be set to the value $1$. It remains to find the
coefficients $\zeta^{(i)}$ and $\la^{(i)}$. By looking at the
coefficient of $w_{i,i+1}$ in $\n_{i+1}(v_i)=r\,(v_i+v_{i+1})$, we
get:
\begin{equation}
r\,\zeta^{(i)}+r^i\la^{(i+1)}=1,\;\;\text{for each $i$ with $1\leq
i\leq n-2$}
\end{equation}
Also, by looking at the coefficient of the same term $w_{i,i+1}$ in
the relation $\n_{i-1}(v_i)=r\,v_i+\unsurr\,v_{i-1}$, we get:
\begin{equation*}
-m\,\zeta^{(i)}-r^{i-3}\,\la^{(i)}=r\,\zeta^{(i)}-\unsur{r^2},\;\;\text{for
each $i$ with $2\leq i\leq n-1$}
\end{equation*}
After multiplication by a factor $r^2$, we obtain:
\begin{equation}
r\,\zeta^{(i)}+r^{i-1}\,\la^{(i)}=1,\;\;\text{for each $i$ with
$2\leq i\leq n-1$}
\end{equation}
By $(22)$ and $(23)$, we get $\la^{(i)}=\unsur{r^{i-2}}\,\la^{(2)}$,
for all $i\geq 2$. Let's do a change of indices in $(22)$ to get:
\begin{equation}r\,\zeta^{(i-1)}+r^{i-1}\la^{(i)}=1\;\;\text{for each $i$ with
$2\leq i\leq n-1$}\end{equation} $(23)$ and $(24)$ show that
$\zeta^{(i)}=\zeta^{(i-1)}$ for each $i$ with $2\leq i\leq n-1$. In
other words, all the $\zeta^{(i)}$ are equal to a certain scalar
$\zeta$. The relation between $\zeta$ and $\la^{(2)}$ is given by
equation $(24)$ with $i=2$:
\begin{equation}\la^{(2)}=\unsurr-\zeta\end{equation} Thus, by
determining $\zeta$, we will get a complete expression for all the
vectors $v_i$'s. Since we have
$$v_1=\zeta\,w_{12}+\sum_{s=3}^n
r^{s-3}\,(w_{1,s}-\unsurr\,w_{2,s}),$$ by looking at the coefficient
of $w_{12}$ in the relation $\n_1(v_1)=-\unsurr\,v_1$, we get the
equation:
\begin{equation}\zeta\bigg(\unsur{l}+\unsurr\bigg)=\unsur{r^2}-(r^2)^{n-3}\end{equation}
Further, by looking at the coefficient of $w_{i,i+1}$ in
$\n_i(v_i)=-\unsurr\,v_i$, we have:
$$\zeta\bigg(\unsur{l}+\unsurr\bigg)=\sum_{s=i+2}^n r^{s-i-3}\,m\,r^{s-i-2}+\la^{(i)}\,\sum_{t=1}^{i-1}r^{t-2}\,\frac{m}{lr^{i-t-1}}$$ \emph{i.e}
$$\zeta\bigg(\unsur{l}+\unsurr\bigg)=\unsur{r^2}-(r^2)^{n-i-2}+\frac{\la^{(i)}}{l}\Big(\unsur{r^i}-r^{i-2}\Big)\qquad(\star)_i$$
\noin Let's write down $(\star)_2$ and $(\star)_3$:
\begin{eqnarray*}
\zeta\,\Big(\unsur{l}+\unsurr\Big)&=&\unsur{r^2}-(r^2)^{n-4}+\frac{\la^{(2)}}{l}\Big(\unsur{r^2}-1\Big)\qquad (\star)_2\\
\zeta\,\Big(\unsur{l}+\unsurr\Big)&=&\unsur{r^2}-(r^2)^{n-5}+\frac{\la^{(2)}}{lr}\Big(\unsur{r^3}-r\Big)\qquad(\star)_3
\end{eqnarray*}
\noin where $\la^{(3)}$ has been replaced by $\frac{\la^{(2)}}{r}$.
Let's subtract these two equalities:
\begin{eqnarray*}
\!\!\!\!\!\!\!\!\frac{\la^{(2)}}{l}\Big(\unsur{r^2}-\unsur{r^4}\Big)\;=\;(r^2)^{n-4}\Big(1-\unsur{r^2}\Big)\qquad
(\star)_2-(\star)_3
\end{eqnarray*}
After multiplying this equality by $\unsur{r^2}$ and dividing it by
$\unsur{r^2}-\unsur{r^4}$ (licit as $m\neq 0$), we obtain:
\begin{equation}
\la^{(2)}=l\,(r^2)^{n-3}
\end{equation}
Hence, by $(25)$, $\zeta=\unsurr-l\,(r^2)^{n-3}$. Plugging this
value for $\zeta$ into $(26)$ now yields:
$$l^2=\unsur{(r^2)^{n-3}},\;\;\text{hence}\;\;l\in\bigg\lb\unsur{r^{n-3}},-\unsur{r^{n-3}}\bigg\rb$$
If $l=\unsur{r^{n-3}}$, we get successively $\la^2=r^{n-3}=\unsur{l}$, $\zeta=\unsurr-\unsur{l}$ and $\la^i=r^{n-i-1}$.\\
If $l=-\unsur{r^{n-3}}$, $\la^2$ and $\zeta$ are still respectively
$\unsur{l}$ and $\unsurr-\unsur{l}$ and $\la^i=-r^{n-i-1}$.
We obtain the formula announced in Theorem $5$. \\
Conversely, if $l\in\lb\unsur{r^{n-3}},-\unsur{r^{n-3}}\rb$, we can
show that the $v_i^{(n)}$'s defined in Theorem $5$ satisfy the
relations $(\tr)$ (see \cite{L}, $\S8.3$). In particular, their
linear span over $F$ is a proper invariant subspace of $\V^{(n)}$,
hence is an $\ih(n)$-module by Corollary $1$. When $n\neq 4$, if the
vectors $v_i^{(n)}$'s were linearly dependent, then their span would
either be one-dimensional or would contain a one-dimensional
$\ih(n)$-submodule, as there is no irreducible $\ih(n)$-module of
dimension between $1$ and $(n-1)$. In any case, by Theorem $4$, that
would force $l=\unsur{r^{2n-3}}$ when $n\neq 3$ and $l\in\lb
-r^3,\unsur{r^3}\rb$ when $n=3$. This is impossible with our
assumption that $l\in\lb\unsur{r^{n-3}}, -\unsur{n-3}\rb$ and the
fact that $r^{2n}\neq 1$. As for $n=4$, the freedom over $F$ of the
family of vectors $(v_1^{(4)},v_2^{(4)},v_3^{(4)})$ is a direct
verification or is a consequence of Theorem $4$ and forthcoming
Proposition $3$ (See $\S3.4$). We are now able to conclude: the
vector space $\text{Span}_F(v_1^{(n)},\dots,v_{n-1}^{(n)})$ is
$(n-1)$-dimensional, is invariant under the action of the $g_i$'s
and is an $\IH_{F,r^2}(n)$-module since it is a proper invariant
subspace of $\V^{(n)}$. Then, by the relations satisfied by
the $v_i^{(n)}$'s, it must be irreducible. \\
To complete the proof of Theorem $5$, we show that there does not
exist any irreducible $5$-dimensional invariant subspace of
$\V^{(6)}$ that is isomorphic to one of the Specht modules
$S^{(3,3)}$ or $S^{(2,2,2)}$. Indeed, suppose such a subspace exists
and name it $\W$. Since we have assumed that $\ih(6)$ is semisimple,
it is licit to use the branching rule as it is described in
Corollary $6.2$ of \cite{M}. We have:
$$S^{(3,3)}\da_{\ih(5)}\,\simeq\, S^{(3,2)}\qquad\qquad
S^{(2,2,2)}\da_{\ih(5)}\,\simeq\, S^{(2,2,1)}$$ We will show that
the restriction of $\W$ to $\ih(5)$ cannot be isomorphic to
$S^{(3,2)}$ or $S^{(2,2,1)}$, hence a contradiction. A proof of the
following fact is in \cite{L}, $\S\,8.3$
\newtheorem{Fact}{Fact}

\begin{Fact}
Suppose $\ih(5)$ is semisimple. Then, up to equivalence, the two
irreducible matrix representations of degree $5$ of $\ih(5)$ are
respectively defined by the matrices $P_1$, $P_2$, $P_3$, $P_4$ and
$Q_1$, $Q_2$, $Q_3$, $Q_4$ given by:
$$P_1:=\begin{bmatrix}
r & & & &\\
  &r& & &\\
  & &r& &\\
1 & &-r^2&-\unsurr&\\
  &1&    &        &-\unsurr\\
\end{bmatrix},\; P_2:=\begin{bmatrix}
-\unsurr& & & 1 & \\
        &-\unsurr&1& & 1\\
        & & r& & \\
        & &  &r&\\
        & &  & &r
\end{bmatrix}$$
$$P_3:=\begin{bmatrix}
r & & & &\\
  &r& & &\\
  &1 &-\unsurr & &\\
1 & & &-\unsurr&-r^2\\
  & & &        &r\\
\end{bmatrix},\; P_4:=\begin{bmatrix}
&1 &-r &  & \\
1 &r-\unsurr&1& & \\
        & & r& & \\
        & & -r^2 & &1\\
        & &  r   & 1 &r-\unsurr
\end{bmatrix}$$
and for the conjugate representation:
$$Q_1:=\begin{bmatrix}
-\unsurr & & & &\\
  &-\unsurr& & &\\
  & &-\unsurr& &\\
1 & &-\unsur{r^2}&r&\\
  &1&    &        &r\\
\end{bmatrix},\; Q_2:=\begin{bmatrix}
r& & & 1 & \\
        &r&1& & 1\\
        & &-\unsurr & & \\
        & &  &-\unsurr&\\
        & &  & &-\unsurr
\end{bmatrix}$$
$$Q_3:=\begin{bmatrix}
-\unsurr & & & &\\
  &-\unsurr& & &\\
  &1 & r & &\\
1 & & &r&-\unsur{r^2}\\
  & & &        &-\unsurr\\
\end{bmatrix},\; Q_4:=\begin{bmatrix}
&1 &\unsurr &  & \\
1 &r-\unsurr&1& & \\
        & &-\unsurr& & \\
        & & -\unsur{r^2} & &1\\
        & &  -\unsurr   & 1 &r-\unsurr
\end{bmatrix}$$

\noin where the blanks must be filled with zeros.
\end{Fact}
First we show that it is impossible to have a basis
$(w_1,w_2,w_3,w_4,w_5)$ of $\W$ in which the matrices of the left
action by the $g_i$'s, $i=1,\dots,4$ are the $Q_i$'s. Indeed,
suppose that such a basis of vectors exists. Let's denote by
$\la_{ij}^{(k)}$ the coefficient of $w_{ij}$ in $w_k$. Since
$g_4.w_4=w_5$ and $g_3.w_5=-\unsur{r^2}\,w_4-\unsurr\,w_5$, we get
$$g_3g_4w_4=-\unsurr\,g_4.w_4-\unsur{r^2}\,w_4$$
We look at the coefficient of $w_{12}$ in this equation to get
$$r^2\,\la_{12}^{(4)}=-\la_{12}^{(4)}-\unsur{r^2}\la_{12}^{(4)}$$
Since $(r^2)^3\neq 0$, we have $r^2+1+\unsur{r^2}\neq 0$ and so
$\la_{12}^{(4)}=0$. Now since $g_3.w_1=-\unsurr\,w_1+w_4$, we get
$r\,\la_{12}^{(1)}=-\unsurr\,\la_{12}^{(1)}$ and so
$\la_{12}^{(1)}=0$. This implies that $\la_{13}^{(1)}$ is also zero
by $g_2.w_1=r\,w_1$ and Lemma $3$. Then, by looking at the
coefficient of $w_{13}$ in $g_2.w_4=w_1-\unsurr\,w_4$, we get
$-m\,\la_{13}^{(4)}=-\unsurr\,\la_{13}^{(4)}$, where we used that
$\la_{12}^{(4)}=0$. Thus, $\la_{13}^{(4)}=0$. Since $g_3.w_4=r\,w_4$
and $g_1.w_4=r\,w_4$, by Lemma $3$, we also get:
$$\la_{13}^{(4)}=\la_{23}^{(4)}=\la_{14}^{(4)}=\la_{24}^{(4)}=0$$
Let's now look at the term in $w_{25}$ in
$g_3g_4\,w_4=-\unsurr\,g_4w_4-\unsur{r^2}\,w_4$. We have:
$$-mr^2\,\la_{15}^{(4)}=-\unsurr\,\bigp-mr\,\la_{15}^{(4)}\bigpd-\unsur{r^2}\,r\,\la_{15}^{(4)},$$
where we used that $\la_{25}^{(4)}=r\,\la_{15}^{(4)}$. Then,
$\la_{15}^{(4)}=0$ and also $\la_{25}^{(4)}=0$. Further, since
$\la_{12}^{(4)}=0$, in $g_1.w_4$, a term in $w_{12}$ is created only
when $g_1$ acts on $w_{26}$, with coefficient $m\,r^3$. Thus the
relation $g_1.w_4=r\,w_4$ yields $\la_{26}^{(4)}=0$. Then by
$g_1.w_4=r\,w_4$, we also have $\la_{16}^{(4)}=0$. \\Furthermore, on
one hand, by looking at the coefficient of $w_{34}$ in
$g_1.w_1=-\unsurr\,w_1+w_4$, we get
$r\,\la_{34}^{(1)}=-\unsurr\,\la_{34}^{(1)}+\la_{34}^{(4)}$,
\emph{i.e\/} $$\la_{34}^{(4)}=(r+\unsurr)\,\la_{34}^{(1)}$$ On the
other hand, by looking at the coefficient of $w_{34}$ in
$g_2.w_4=w_1-\unsurr\,w_4$ and remembering that $\la_{24}^{(4)}=0$,
we have
$-m\,\la_{34}^{(4)}=-\unsurr\,\la_{34}^{(4)}+\la_{34}^{(1)}$,
\emph{i.e\/} $$\la_{34}^{(4)}=\unsurr\,\la_{34}^{(1)}$$ The two
relations binding $\la_{34}^{(4)}$ and $\la_{34}^{(1)}$ now yield
$\la_{34}^{(4)}=\la_{34}^{(1)}=0$. \\
So $w_4$ reduces to
$$w_4=r\la_{35}^{(4)}\,w_{45}+\la_{35}^{(4)}\,w_{35}+
\la_{56}^{(4)}\,w_{56}+ r\la_{36}^{(4)}\,w_{46}+
\la_{36}^{(4)}\,w_{36}$$ \noin At this point, it is of interest to
derive the following result.
\newtheorem{Result}{Result}
\begin{Result}
The irreducible matrix representation of degree $5$ of $\ih(5)$
defined by the matrices $Q_i$'s is not a constituent of the
Lawrence-Krammer representation of degree $10$ of the BMW algebra of
type $A_4$.
\end{Result}
\noin\textbf{Proof.} If $\W$ is a subspace of $\V^{(5)}$ instead,
then we simply have
$$w_4=r\la_{35}^{(4)}\,w_{45}+\la_{35}^{(4)}\,w_{35}$$
Then, by looking at the coefficient of $w_{34}$ in $g_3.w_4=r\,w_4$,
we get $m\,\la_{45}^{(4)}=0$, so that
$\la_{35}^{(4)}=\la_{45}^{(4)}$. Then $w_4=0$, which is impossible.
$\;\;\;\square$\\\\
Let's go back to the main proof. By looking at the coefficient of
$w_{34}$ in $g_3.w_4=r\,w_4$ and using the complete expression for
$w_4$ this time, we get: $m\,r\la_{35}^{(4)}+m\,r^2\la_{36}^{(4)}$,
\emph{i.e\/} $\la_{36}^{(4)}=-\unsurr\,\la_{35}^{(4)}$. We see that
all the coefficients in $w_4$ except $\la_{56}^{(4)}$ are multiples
of $\la_{35}^{(4)}$. Moreover, we claim that $w_4$ may not be a
mutiple of $w_{56}$. Indeed, recall the formula
$$g_3g_4.w_4=-\unsurr\,g_4.w_4-\unsur{r^2}\,w_4$$ and observe that in
$-\unsurr\,g_4.w_4-\unsur{r^2}\,w_4$, there is no term in $w_{36}$
while in $g_3g_4.w_4$, there is one. Thus, without loss of
generality, we may set $\la_{35}^{(4)}=1$. So,
$$w_4=r\,w_{45}+w_{35}+\la_{56}^{(4)}\,w_{56}-w_{46}-\unsurr\,w_{36}$$
We then deduce a complete expression for $w_1$ by using the relation
$w_1=g_2.w_4+\unsurr\,w_4$:
$$w_1=(1+r^2)\,w_{45}+r\,w_{35}+w_{25}\\+\bigp
r+\unsurr\bigpd\,\la_{56}^{(4)}\,w_{56}-\bigp
r+\unsurr\bigpd\,w_{46}-w_{36}-\unsurr\,w_{26}$$ A contradiction now
arises when looking at the coefficient of $w_{26}$ in
$g_3.w_1=-\unsurr\,w_1+w_4$. Indeed, this yields $-1=\unsur{r^2}$
and contradicts $(r^2)^2\neq 1$.\\

Suppose now that there exists a basis $(w_1,w_2,w_3,w_4,w_5)$ of
$\W$ in which the matrices of the left action by the $g_i$'s,
$i=1,\dots,4$ are the $P_i$'s. We read on the matrices $P_1$ and
$P_3$ that $g_1.w_4=-\unsurr\,w_4$ and $g_3.w_4=-\unsurr\,w_4$.
Thus, we have:
$$w_4=\mu_{13}^{(4)}w_{13}+\mu_{23}^{(4)}\,w_{23}+\mu_{24}^{(4)}\,w_{24}+\mu_{14}^{(4)}\,w_{14},$$
where the coefficients are related by
$\mu_{14}^{(4)}=\mu_{23}^{(4)}=-\unsurr\,\mu_{13}^{(4)}=-r\,\mu_{24}^{(4)}$.
In particular, all these coefficients are nonzero. As
$g_1.w_1=g_3\,w_1$, by looking at the coefficient in $w_{23}$, we
obtain
$$\mu_{13}^{(1)}-m\,\mu_{23}^{(1)}=\mu_{24}^{(1)}$$ Moreover, by
looking at the same coefficient $w_{23}$ in $g_3.w_1=r\,w_1+w_4$, we
get
$$\mu_{24}^{(1)}=r\,\mu_{23}^{(1)}+\mu_{23}^{(4)}$$ Combining both
equations yields
$$\mu_{13}^{(1)}-\unsurr\,\mu_{23}^{(1)}=\mu_{23}^{(4)}$$ Further, by
looking at the coefficient in $w_{13}$ in the same equation, we get
$0=r\mu_{13}^{(1)}+\mu_{13}^{(4)}$, as there is no term in $w_{14}$
in $w_1$ by $g_2.w_1=-\unsurr\,w_1$. Then we get $\mu_{23}^{(1)}=0$.
Looking at the coefficient of $w_{23}$ in $g_2.w_4=w_1+r\,w_4$ now
yields
$r\,\mu_{23}^{(4)}=\unsur{l}\,\mu_{23}^{(4)}+\frac{m}{l}\,\mu_{13}^{(4)}$,
which by using the relations from the beginning reads
$$\bigg(\frac{m}{l}+\unsurr(1-\unsur{l}\bigg)\,\mu_{13}^{(4)}=0$$
As the coefficient $\mu_{13}^{(4)}$ is nonzero, this forces $l=r$.
It will be useful to derive the following result along the way:
\begin{Result}
Assume $\ih(5)$ is semisimple. If there exists an irreducible
$5$-dimensional invariant subspace of $\V^{(5)}$ then $l=r$.
\end{Result}
\noin\textbf{Proof.} Indeed, if we assume $\W\subset\V^{(5)}$
instead of $\W\subset\V^{(6)}$ in the computations above, they are
unchanged and lead to the same conclusion. $\;\;\;\square$\\\\
As seen along the way, $w_4$ is, up to a multiplication by a scalar
that can be set to $1$ without loss of generality,
$$w_4=w_{23}+w_{14}-r\,w_{13}-\unsurr\,w_{24}$$ Then the other
spanning vectors of $\W$ must be:
$$\left.\begin{array}{ccccc}
w_5&=&g_4.w_4&=&r\,w_{23}+w_{15}-r^2\,w_{13}-\unsurr\,w_{25}\\
w_1&=&g_2.w_4-r\,w_4&=&w_{24}-r\,w_{12}+w_{13}-\unsurr\,w_{34}\\
w_2&=&g_4.w_1&=&w_{25}-r^2\,w_{12}+r\,w_{13}-\unsurr\,w_{35}\\
w_3&=&g_3.w_2-r\,w_2&=&r\,w_{14}-r^2\,w_{13}+w_{35}-\unsurr\,w_{45}
\end{array}\right.$$

\noin where we replaced $l$ by $r$. But $\W$ is an invariant
subspace of $\V^{(6)}$. In particular, it must be invariant under
the action by $g_5$. This is not compatible with the spanning set
above. We conclude that it is impossible to have
$$\W\da_{\ih(5)}\,\simeq\, S^{(3,2)}\qquad\text{or}\qquad
\W\da_{\ih(5)}\,\simeq\, S^{(2,2,1)}$$ and so $\W$ cannot be
isomorphic to $S^{(3,3)}$ or $S^{(2,2,2)}$. Thus, by previous work,
the existence of an irreducible $5$-dimensional invariant subspace
of $\V^{(6)}$ implies that $l\in\lb\unsur{r^3},-\unsur{r^3}\rb$.
This completes the proof of the Theorem. $\;\;\;\square$

\subsection{The cases $l=r$ and $l=-r^3$}
In this section, we show that when $l=r$ the representation
$\n^{(n)}$ is reducible for all $n\geq 4$ and when $l=-r^3$, the
representation is reducible for all $n\geq 3$. The latter point is
true by Theorem $4$ when $n=3$ and by Theorem $5$ when $n=4$. To do
so, we show that some proper invariant subspace $K(n)$ of
$\V^{(n)}$, defined in the Proposition below, is nontrivial.
\begin{Proposition}
For any two nodes $i$ and $j$ with $1\leq i<j\leq n$, define
$$c_{ij}=\begin{cases}g_{j-1}\dots g_{i+1}e_ig_{i+1}^{-1}\dots g_{j-1}^{-1}&\text{if $j\geq i+2$}\\
e_i&\text{if $j=i+1$}\end{cases}$$ Then, $K(n)=\cap_{1\leq i<j\leq
n}\text{Ker}\;\n^{(n)}(c_{ij})$ is a proper invariant subspace of
$\V^{(n)}$. Moreover, any proper invariant subspace of $\V^{(n)}$
must be contained in $K(n)$.
\end{Proposition}
\noin\textbf{Proof.} $K(n)$ is proper, as is visible on the
expressions for $\n^{(n)}(e_i)$. Further, if an $\xb$ is annihilated
by all the $g_i$ conjugates of the $e_i$'s, then $\n_k(\xb)$ is also
annihilated by these same elements. Verification of this fact is
tedious and can be found in \cite{L}, $\S2$. Hence $K(n)$ is
invariant. Let $\W$ be a proper invariant subspace of $\V^{(n)}$. By
Proposition $1$, we have $\n^{(n)}(c_{i,i+1})(\W)=0$ for all $i$
with $1\leq i\leq n-1$. This fact is also true for the other
conjugates $c_{ij}$'s. Hence $\W$ must be contained in $K(n)$
$\;\;\;\square$
\\\\
This is how an element $c_{ij}$ is represented in the tangle
algebra:\vspace{-0.1in}
\begin{center}
\epsfig{file=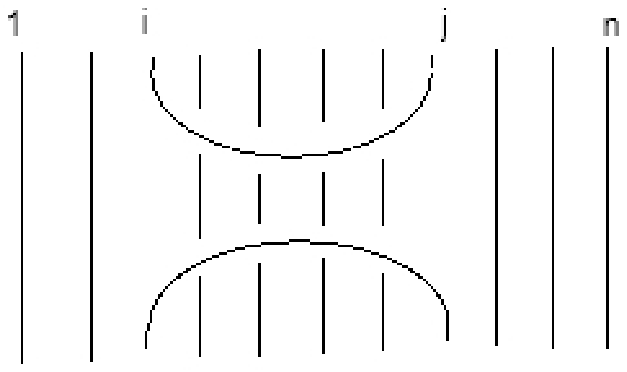, height=3.5cm}\end{center}\vspace{-0.25in} It
has two horizontal strands: one at the top and one at the bottom,
each joining nodes $i$ and $j$ and moreover, when $j\geq i+2$, such
horizontal strands over-cross all the vertical strands that they
intersect.

To show that $\n^{(n)}$ is reducible, it will suffice to exhibit a
nontrivial element in $K(n)$ when $l=r$ or $l=-r^3$. The following
Proposition shows that $K(4)$ is nontrivial and irreducible when
$l=r$ and $\ih(4)$ is semisimple.
\begin{Proposition}
Assume $\ih(4)$ is semisimple. There exists an irreducible
$2$-dimensional invariant subspace of $\V^{(4)}$ if and only if
$l=r$. If so it is unique and it is $K(4)$. Moreover, it is spanned
over $F$ by the two linearly independent vectors:
\begin{eqnarray}
v_1&=&w_{13}-\unsurr\,w_{23}+\unsur{r^2}\,w_{24}-\unsurr\,w_{14}\\
v_2&=&w_{12}-\unsurr\,w_{13}-\unsurr\,w_{24}+\unsur{r^2}\,w_{34}
\end{eqnarray}
\end{Proposition}
\noin\textbf{Proof.} When $\ih(4)$ is semisimple, the following
three matrices
$$ H_1=\begin{bmatrix}
-\unsurr & 1\\
0 & r\\
\end{bmatrix},\; H_2=\begin{bmatrix}
r &0\\
1&-\unsurr\\
\end{bmatrix},\; H_3=\begin{bmatrix}
-\unsurr & 1\\
0 & r\\ \end{bmatrix}$$ define an irreducible matrix representation
of degree $2$ of $\ih(4)$. Suppose $\W$ is an irreducible
$2$-dimensional invariant subspace of $\V^{(4)}$. Then $\W$ has a
basis $(v_1,v_2)$ of vectors such that the matrix of $\n_i$ in this
basis is $H_i$. Since $\n_1(v_1)=-\unsurr\,v_1$ (resp
$\n_3(v_1)=-\unsurr\,v_1$), there is no term in $w_{34}$ (resp
$w_{12}$) in $v_1$ and since $\n_2(v_2)=-\unsurr\,v_2$, there is no
term in $w_{14}$ in $v_2$. Let the $\la_{ij}$'s (resp $\mu_{ij}$'s)
denote the coefficients of the $w_{ij}$'s in $v_1$ (resp $v_2$). We
have:
$$\left\lb\begin{array}{cccc}\la_{23}=-\unsurr\,\la_{13}&\text{and}&\la_{24}=-\unsurr\,\la_{14}&\text{by $(14)$ with $i=1$ and $\gamma=-\unsurr$}\\
\la_{24}=-\unsurr\,\la_{23}&\text{and}&\la_{14}=-\unsurr\,\la_{13}&\text{by
$(15)$ with $i=3$ and $\gamma=-\unsurr$}\\
\mu_{34}=-\unsurr\,\mu_{24}&&&\text{by $(14)$ with $i=2$ and
$\gamma=-\unsurr$}\\
\mu_{13}=-\unsurr\,\mu_{12}&&&\text{by $(15)$ with $i=2$ and
$\gamma=-\unsurr$}\end{array}\right.$$ Hence, without loss of
generality,
$v_1=w_{13}-\unsurr\,w_{23}+\unsur{r^2}\,w_{24}-\unsurr\,w_{14}$ and
$v_2$ is a multiple of
$w_{12}-\unsurr\,w_{13}+\mu(w_{24}-\unsurr\,w_{34})+\mu^{'}\,w_{23},$
where $\mu$ and $\mu^{'}$ are scalars to determine. The relation
$\n_2(v_1)=r\,v_1+v_2$ sets
$v_2=w_{12}-\unsurr\,w_{13}+\mu(w_{24}-\unsurr\,w_{34})+\mu^{'}\,w_{23},$
by just looking at the coefficient in $w_{12}$. The same relation
yields $\mu=-\unsurr$ by looking at the coefficient in $w_{24}$.
Next, by looking at the coefficient of $w_{23}$ in
$\n_3(v_2)=v_1+r\,v_2$, we get $r\,\mu^{'}-\unsurr=\mu$. Replacing
$\mu$ by $-\unsurr$ now yields $\mu^{'}=0$. We thus get the
expressions in $(28)$ and $(29)$ for $v_1$ and $v_2$ respectively.
Also, by looking at the coefficient in $w_{12}$ in
$\n_1(v_2)=v_1+r\,v_2$, we have $\unsur{l}-m=r$, hence $l=r$.
Conversely, if $l=r$, it is a direct verification that the vectors
$v_1$ and $v_2$ given by the formulas $(28)$ and $(29)$ are linearly
independent and that they verify the relations:
\begin{eqnarray*}
\n_i(v_i)&=&-\unsurr\,v_i\;\;\;\text{when}\;\; i\in\lb 1,2\rb
\;\;\;\&\;\;\;\n_3(v_1)\;\;=\;\;-\unsurr\,v_1\\
\n_1(v_2)&=&v_1+r\,v_2\;\;=\;\;\n_3(v_2)\;\;\&\;\;\n_2(v_1)\;\;
=\;\; r\,v_1+v_2
\end{eqnarray*}
Thus, their linear span over $F$ is an irreducible $2$-dimensional
invariant subspace of $\V^{(4)}$. It remains to show that it is in
fact $K(4)$. By Proposition $2$, $Span_F(v_1,v_2)$ is contained in
$K(4)$. If $K(4)$ is three-dimensional, either it is irreducible and
so $l\in\lb -r^3,\unsurr,-\unsurr\rb$ by Theorem $5$. This is
impossible as $l=r$. Or it is reducible and it must contain a
one-dimensional invariant subspace. Then by Theorem $4$, it forces
$l=\unsur{r^5}$, which is again impossible. If $K(4)$ is
four-dimensional, then $K(4)$ is not irreducible as its dimension is
not $1,2$ or $3$. Since we just saw that there exists only one
irreducible $2$-dimensional invariant subspace of $\V^{(4)}$ when
$l=r$, $K(4)$ must then contain a one-dimensional invariant
subspace, which is again impossible. For similar reasons, it is also
impossible to have $k(4)=5$, hence the only possibility that is left
is to have $k(4)=2$. We conclude that $K(4)=Span_F(v_1,v_2)$ and
$K(4)$ is thus
irreducible. $\;\;\;\square$\\\\
The next Proposition shows the reducibility of the representation
when $l=r$ and $n\geq 4$, where we still assume that $\ih(n)$ is
semisimple.
\begin{Proposition}
Assume $l=r$ and $\ih(n)$ is semisimple. \\Then the vector
$v_1=w_{13}-\unsurr\,w_{23}+\unsur{r^2}\,w_{24}-\unsurr\,w_{14}$ of
Proposition $3$ belongs to $K(n)$ for all $n\geq 4$. Thus,
$\n^{(n)}$ is reducible when $l=r$ and $n\geq 4$.
\end{Proposition}
\noin\textbf{Proof.} For $n=4$, the result is contained in
Proposition $3$. When $i\geq 5$, we simply have for any $j\geq i+2$,
$\n_{i+1}^{-1}\dots\n_{j-1}^{-1}(v_1)=\unsur{r^{j-i-1}}\,v_1$ and
since $\n^{(n)}(e_i)(v_1)=0$, we see that $v_1$ is thus annihilated
by all the $\n^{(n)}(c_{ij})$ with $i\geq 5$. Also, since we just
saw in Proposition $3$ that $v_1$ is in $K(4)$, $v_1$ is annihilated
by all the $\n^{(n)}(c_{ij})$'s with $j\leq 4$. Thus, it suffices to
check that $v_1$ is annihilated by $\n^{(n)}(c_{1j})$,
$\n^{(n)}(c_{2j})$, $\n^{(n)}(c_{3j})$ and $\n^{(n)}(c_{4j})$ for
any $j\geq 5$. We will use the following formulas that give the
action of the $c_{ij}$'s on the basis vectors of the L-K space in
some relevant cases here:
$$\left|\begin{array}{l}
\n^{(n)}(c_{ij})(w_{i,j-k})=\unsur{lr^{k-1}}\,w_{ij}\qquad\qquad\qquad\qquad\qquad\;\!\!\!(R_k)_{1\leq k\leq j-i-1}\\
\n^{(n)}(c_{ij})(w_{i-k,i})=\unsur{r^{(k-1)+(j-i-1)}}\,w_{ij}\qquad\qquad\qquad\;\;\;\!\!(L_{j-i,k})_{1\leq k\leq i-1}\\
\n^{(n)}(c_{ij})(w_{i-t,j-s})=(\unsur{r^{t+s-1}}-\unsur{r^{t+s-3}})(\unsur{l}-\unsurr)\,w_{ij}\!\!\;\;\;\;\,(C_{t,s})_{1\leq
t\leq i-1,\;1\leq s\leq j-i-1}
\end{array}\right.$$

\noin These formulas can be shown and pictured easily by using the
tangles. Let's take an example in $\V^{(12)}$. The product tangle
$c_{4,9}\,w_{2,7}$ as represented in the figure below
\vspace{-0.2in}
\begin{center} \epsfig{file=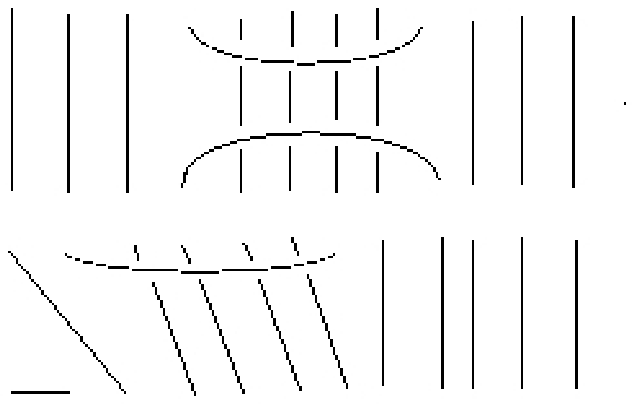,
height=4cm}\end{center}\vspace{-0.03in} expands as follows, where we
use the Kauffman skein relation:\vspace{-0.03in}
$$\begin{array}{l}
\text{\epsfig{file=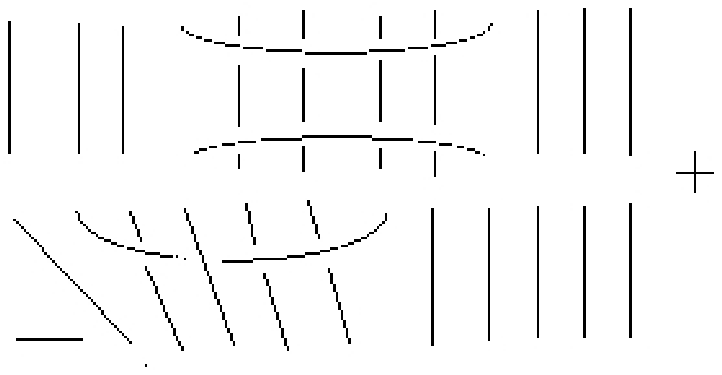, height=4cm}}\text{\epsfig{file=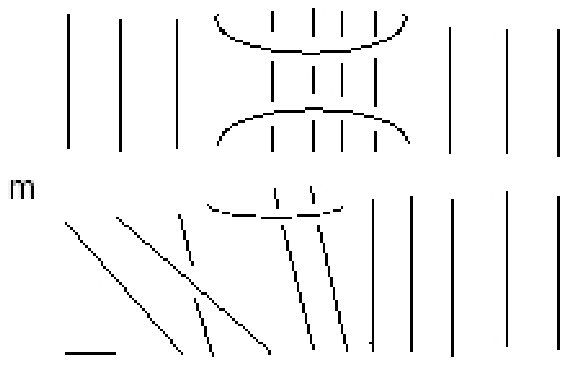,
height=4cm}}\\
\qquad\qquad\qquad\qquad\text{\epsfig{file=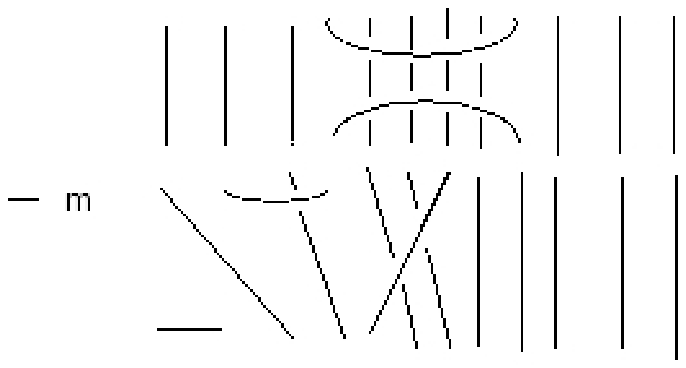, height=4cm}}
\end{array}\vspace{-0.25in}$$
After doing a Reidemeister's move of type \setcounter{c}{2}
\Roman{c} (as it is described in $\S2.2$ of \cite{MW}) on the first
tangle of the sum above, we see that it is zero. After "delooping"
the second tangle of the sum and using Reidemeister's move \Roman{c}
twice, we see that it is obtained from the basis vector $w_{4,9}$
with non-crossed vertical strands by multiplying it to the right by
$g_4^{-1}g_8^{-1}$. The resulting coefficient is
$\frac{m}{lr^2}=\unsur{l}(\unsur{r^3}-\unsurr)$. If we call
over-crossing a multiplication at the bottom by a $g_i$ in $H$ and
under-crossing a multiplication at the bottom by a $g_i^{-1}$ in
$H$, we see that in the last tangle of the sum, there are five
under-crossings and two over-crossings. Thus, the resulting
coefficient is $-\frac{mr^2}{r^5}=\unsur{r^2}-\unsur{r^4}$. When
adding these two coefficients, we get
$(\unsur{r^3}-\unsurr)(\unsur{l}-\unsurr)$, which the coefficient in
$(C_{2,2})$. Using this example as a support, it is easy to see
that, more generally, for any fixed nodes $i$ and $j$ with $1\leq
i<j\leq n$, the coefficient of $\n^{(n)}(c_{ij})(w_{i-t,j-s})$,
where $1\leq t\leq i-1$ and $1\leq s\leq j-i-1$, is:
$$m\,\bigg(\unsur{l}\unsur{r^{(t-1)+(s-1)}}-\frac{r^{\lb(j-s)-i\rb-1}}{r^{(t-1)+(j-i-1)}}\bigg)$$
which after simplification is the coefficient in $(C_{t,s})$. Thus,
we obtain the family of equations $(C_{t,s})_{1\leq t\leq
i-1,\,1\leq s\leq j-i-1}$. Similarly for fixed nodes $i$ and $j$,
the two families of equations $(R_k)_{1\leq k\leq j-i-1}$ and
$(L_{j-i,k})_{1\leq k\leq i-1}$ can be pictured easily by using the
tangles, or can be established by using the definition of $c_{ij}$
and the expression for the representation in $\S2.3.2$.\\
When $l=r$, we note that the action of $c_{i,j}$ on $w_{i-t,j-s}$ is
zero. From there, we have for $j\geq 5$, where we replaced $l$ by
$r$:
$$\left\lb\begin{array}{l}
\n^{(n)}(c_{1,j})(v_1)=\n^{(n)}(c_{1,j})(w_{13}-\unsurr\,w_{14})=0\;\;\qquad\text{by
$(R_{j-3})$ and $(R_{j-4})$}\\
\n^{(n)}(c_{2,j})(v_1)=\n^{(n)}(c_{2,j})(-\unsurr\,w_{23}+\unsur{r^2}\,w_{24})=0\;\;\,\text{by
$(R_{j-3})$ and $(R_{j-4})$}\\
\n^{(n)}(c_{3,j})(v_1)=\n^{(n)}(c_{3,j})(w_{13}-\unsurr\,w_{23})=0\;\;\qquad\text{by
$(L_{j-3,2})$ and $(L_{j-3,1})$}\\
\n^{(n)}(c_{4,j})(v_1)=\n^{(n)}(c_{4,j})(\unsur{r^2}\,w_{24}-\unsurr\,w_{14})=0\;\;\;\;\;\text{by
$(L_{j-4,2})$ and $(L_{j-4,3})$}
\end{array}\right.$$
\noin So $v_1$ is in $K(n)$ for all $n\geq 4$, as announced. It will
be useful to note on the way that by the game of the coefficients,
the equalities to the right of the first two lines of equations
still hold when $l=-r^3$.
$\;\;\;\square$\\
When $l=-r^3$, we have a similar result. This is the object of the
next proposition.
\begin{Proposition}
Assume $\ih(n)$ is semisimple. When $l=-r^3$, the vector $u_1$
defined as in Theorem $5$ by
$u_1=r\,w_{23}+w_{13}+(\unsurr+\unsur{r^3})w_{34}-w_{24}-\unsurr\,w_{14}$
belongs to $K(n)$ for all $n\geq 4$. Thus, when $l=-r^3$, the
representation $\n^{(n)}$ is reducible for every $n\geq 3$.
\end{Proposition}
\noin\textbf{Proof.} When $l=-r^3$, $\n^{(3)}$ is reducible by
Theorem $4$ and $\n^{(4)}$ is also reducible by Theorem $5$. Suppose
now $n\geq 5$. To show that $u_1$ is in $K(n)$, like in the case
$l=r$, it will suffice to check that $\n^{(n)}(c_{ij})(u_1)=0$ for
all $i\leq 4$ and $j\geq 5$. With $l=-r^3$, the coefficients of type
$(C_{t,s})$ are no longer zero. But we have:
$\n^{(n)}(c_{2,j})(w_{13}-\unsurr\,w_{14})=0$ by $(C_{1,j-3})$ and
$(C_{1,j-4})$. For $\n^{(n)}(c_{3,j})(u_1)$, there is no shortcut
and a complete evaluation must be performed. We have, where we
respected the same order of the terms in Proposition $5$ for the
coefficients:
\begin{equation*}\begin{split}\n^{(n)}(c_{3,j})&(u_1)=r\,\unsur{r^{j-4}}+\unsur{r^{j-3}}+\bigg(\unsurr+\unsur{r^3}\bigg)\bigg(-\unsur{r^3\,r^{j-5}}\bigg)\\&+\bigg(\unsur{r^{j-4}}-\unsur{r^{j-6}}\bigg)\bigg(\unsur{r^3}+\unsurr\bigg)+\unsurr\bigg(\unsur{r^{j-3}}-\unsur{r^{j-5}}\bigg)\bigg(\unsur{r^3}+\unsurr\bigg)\,w_{3,j}\end{split}\end{equation*}
The rules used are, in the same order: $(L_{j-3,1})$, $(L_{j-3,2})$,
$(R_{j-4})$, $(C_{1,j-4})$ and $(C_{2,j-4})$. All the coefficients
cancel nicely to give $\n^{(n)}(c_{3,j})(u_1)=0$.\\
Finally, for $\n^{(n)}(c_{4,j})(u_1)$, only the terms in $u_1$ whose
last node is node number $4$ yield a nonzero contribution, the first
one contributing with a coefficient
$(\unsurr+\unsur{r^3})\,\unsur{r^{j-5}}$, the second one with a
coefficient $-\unsur{r^{j-4}}$ and the third one with a coefficient
$-\unsurr\,\unsur{r^{j-3}}$ by rules $(L_{j-4,1})$, $(L_{j-4,2})$
and $(L_{j-4,3})$ respectively. The sum of these three coefficients
is zero. \\
Thus, we are done with all the cases and $u_1$ belongs to $K(n)$ for
all $n\geq 4$. $\;\;\;\square$\\
At this stage, we have shown that when $l$ and $r$ take the values
of Theorem $1$, the representation $\n^{(n)}$ is reducible. In the
next section, we show conversely that if $\n^{(n)}$ is reducible,
then $l$ and $r$ must related in the way described in Theorem $1$.
\section{Proof of the main theorem}
We recall from Proposition $2$ that any proper irreducible invariant
subspace of $\V^{(n)}$ is an irreducible $\ih(n)$-module. When
$n=3$, the irreducible $\ih(3)$-modules have dimension $1$ or $2$.
We showed in Theorem $4$ that there exists a one-dimensional
invariant subspace of $\V^{(3)}$ if and only if
$l\in\lb\-r^3,\unsur{r^3}\rb$ and we saw in Theorem $5$ that there
exists an irreducible $2$-dimensional invariant subspace of
$\V^{(3)}$ if and only if $l\in\lb 1,-1\rb$. Hence the main theorem
is proven for $n=3$. When $n=4$, the irreducible $\ih(4)$-modules
have dimensions $1$, $2$ or $3$. By Theorem $4$ (resp Theorem $5$,
resp Proposition $3$), there exists a one-dimensional (resp an
irreducible $2$-dimensional, resp an irreducible $3$-dimensional)
invariant subspace of $\V^{(4)}$ if and only if $l=\unsur{r^5}$
(resp $l=r$, resp $l\in\lb -r^3, \unsurr, -\unsurr\rb$). Thus, the
main theorem also holds for $n=4$. Suppose now $n\geq 5$. By $\S3$,
it suffices to prove that if $\n^{(n)}$ is reducible, then $l\in\lb
r,-r^3, \unsur{r^{2n-3}}, \unsur{r^{n-3}}, -\unsur{r^{n-3}}\rb$ and
the proof of the main theorem will be complete. We show it for $n=5$
and $n=6$, then proceed by induction on $n$.
\subsection {The case $n=5$}
If $\W$ is an irreducible proper invariant subspace of $\V^{(5)}$,
then $\di(\W)\in\lb 1,4,5,6\rb$, as $\W$ is an irreducible
$\ih(5)$-module by Corollary $1$ of $\S3.1$. If $\di(\W)=5$, it
forces $l=r$ by Result $2$. If $\di(\W)=1$, it forces
$l=\unsur{r^7}$ by Theorem $4$ and if $\di(\W)=4$, it forces
$l\in\lb\unsur{r^2},-\unsur{r^2}\rb$ by Theorem $5$. From now on,
assume that $l\not\in\lb r,\unsur{r^7},\unsur{r^2},-\unsur{r^2}\rb$.
We will show that $l=-r^3$. By choice for $l$ and $r$, we have
$\di(\W)=6$. Then $\W\cap\V^{(4)}\neq \lb 0\rb$, as otherwise
$\W\oplus\V^{(4)}\subseteq\V^{(5)}$, which implies on the dimensions
$\di\,\W\leq 10-6=4$. Since $\W\cap\V^{(4)}$ is a proper invariant
subspace of $\V^{(4)}$, the representation $\n^{(4)}$ is then
reducible, which yields $l\in\lb
-r^3,\unsurr,-\unsurr,\unsur{r^5}\rb$. We will show that it is
impossible to have $l\in\lb\unsurr,-\unsurr,\unsur{r^5}\rb$, unless
$l=\unsur{r^5}=-r^3$.

Let's first assume that $l=\unsur{r^5}$ and show that under our
assumptions, it forces $l=-r^3$. When $l=\unsur{r^5}$, there exists
a one-dimensional invariant subspace, say $\U_1$, of $\V^{(4)}$ by
Theorem $4$ and by Proposition $2$ it is contained in $K(4)$. In
particular the dimension $k(4)$ of $K(4)$ is $1,2,3,4$ or $5$. But
since $\W\cap\V^{(4)}$ is a proper invariant subspace of $\V^{(4)}$,
it must be contained in $K(4)$ by Proposition $2$. Hence we have
$k(4)\geq \di(\W\cap\V^{(4)})$. Also, since
$\di(\W\cap\V^{(4)})=\di(\W)+\di(\V^{(4)})-\di(\W+\V^{(4)})\geq
12-\di(\V^{(5)})=2$, we get $k(4)\geq 2$. By semisimplicity of
$\ih(4)$ and by uniqueness of a one-dimensional invariant subspace
of $\V^{(4)}$ when it exists, $\U_1$ then has a summand $S$ in
$K(4)$ that is of dimension greater than or equal to $2$. $S$ is an
invariant subspace of $\V^{(4)}$ and it cannot contain a
one-dimensional invariant subspace. Nor can it contain an
irreducible $2$-dimensional invariant subspace by Proposition $3$
since we assumed $l\neq r$. Also, since $1+\di(S)=k(4)\leq 5$, we
note that $\di(S)\leq 4$. Then, by the same arguments already
exposed, $\di(S)=3$ and $S$ is irreducible. By Theorem $5$, we now
get $l\in\lb -\unsurr,\unsurr,-r^3\rb$. But since $l=\unsur{r^5}$,
it forces $l=-r^3$ as it is impossible to have $(r^2)^4=1$ when
$\ih(5)$ is semisimple.

Assume next that $l\in\lb\unsurr,-\unsurr\rb$. We show that these
values lead to a contradiction. First, by choice for $l$ and $r$ and
Theorem $5$ (case $n=4$), $\V^{(4)}$ contains an irreducible
$3$-dimensional invariant subspace and by Proposition $2$, this
proper invariant subspace must be contained in $K(4)$. Hence
$k(4)\geq 3$. Since there cannot exist any one-dimensional invariant
subspace of $\V^{(4)}$ (as $l\neq\unsur{r^5}$ when $(r^2)^4\neq 1$)
or any irreducible $2$-dimensional invariant subspace of $\V^{(4)}$
(as $l\neq r$), we cannot have $k(4)\in\lb 4,5\rb$. Thus, we have
$k(4)=3$ and so $K(4)$ is irreducible. Now the irreducibility of
$K(4)$ and the fact that $0\subset \W\cap\V^{(4)}\subseteq K(4)$
implies that
\begin{equation}\W\cap\V^{(4)}=K(4)\end{equation}
We show that this is impossible.\\
Consider first the case when $l=\unsurr$. When $l=\unsurr$, the
vector $w_{14}-w_{23}$ belongs to $K(4)$. Indeed, this vector is
$$\unsur{r^2+\unsur{r^2}}\bigg(r\,v_1^{(4)}+(r-\unsurr)\,v_2^{(4)}-\unsurr\,v_3^{(4)}\bigg)$$
Then it also belongs to $\W$. It follows that
$\n^{(5)}(e_4)(w_{14}-w_{23})=\unsur{r^2}\,x_{\al_4}\in\W$, as $\W$
is an invariant subspace of $\V^{(5)}$. But then, by the same
argument as in the proof of Proposition $1$, $\W$ is the whole space
$\V^{(5)}$, in contradiction with $\W$ is proper. Thus, it is
impossible to have $(30)$ and so $l$ cannot take the value
$\unsurr$.\\
Consider now the case when $l=-\unsurr$. The vector
$w_{14}-w_{23}+w_{12}+w_{34}$ belongs to $K(4)$ since it is
$$\unsur{r+\unsurr}\bigg(v_1^{(4)}+v_3^{(4)}\bigg)$$
By $(30)$, this vector also belongs to $\W$. But then
$$\n^{(5)}(e_4)(w_{14}-w_{23}+w_{12}+w_{34})=\bigg(1+\unsur{r^2}\bigg)\,x_{\al_4}\in\W$$
As $(r^2)^2\neq 1$, this implies in turn that $x_{\al_4}$ is in
$\W$. Then $\W$ is the whole space $\V^{(5)}$, a contradiction.
Thus, it is also impossible to have $l=-\unsurr$.

We have now shown that if $\n^{(5)}$ is reducible and $l\not\in\lb
r, \unsur{r^7},\unsur{r^2},-\unsur{r^2}\rb$, then $l=-r^3$. Thus, if
$\n^{(5)}$ is reducible, then $l\in\lb r,-r^3,\unsur{r^7},
-\unsur{r^2},\unsur{r^2}\rb$.
\subsection{The case $n=6$}
Let $\W$ be an irreducible proper invariant subspace of $\V^{(6)}$.
So $\W$ is an irreducible $\ih(6)$-module. The irreducible
representations of $\ih(6)$ have degrees $1$, $5$, $9$, $10$, $16$.
The vector space $\V^{(6)}$ is $15$-dimensional. Hence
$\di(\W)\in\lb 1,5,9,10\rb$. If $\di(\W)=1$, then $l=\unsur{r^9}$ by
Theorem $4$ and if $\di(\W)=5$, then
$l\in\lb\unsur{r^3},-\unsur{r^3}\rb$ by Theorem $5$. Suppose now
$l\not\in\lb \unsur{r^9}, \unsur{r^3}, -\unsur{r^3}\rb$. Then
$\di(\W)\geq 9$, which implies in particular that
$\W\cap\V^{(5)}\neq \lb 0\rb$. Moreover, $\W\cap\V^{(5)}$ is a
proper subspace of $\V^{(5)}$, as otherwise $\W$ would contain
$\V^{(5)}$ and would in fact be the whole space $\V^{(6)}$. Hence we
see that $\n^{(5)}$ is reducible, which implies that
\begin{equation}l\in\bigg\lb
r,-r^3,\unsur{r^7},\unsur{r^2},-\unsur{r^2}\bigg\rb\end{equation} by
the case $n=5$. Also, if $\W\cap\V^{(4)}=\lb 0\rb$, then
$\W\oplus\V^{(4)}\subseteq\V^{(6)}$ and so $\di(\W)\leq 15-6=9$.
Then $\di\W=9$. We notice that
$\di(\W)+\di(\V^{(4)})=\di(\V^{(6)})$. Thus, we get
$\W\oplus\V^{(4)}=\V^{(6)}$. But since $\W\subseteq K(6)$ by
Proposition $2$, we must have in particular $\n^{(6)}(e_5)(\W)=0$.
But $e_5$ also acts trivially on $\V^{(4)}$. It follows that $e_5$
acts trivially on the direct sum $\W\oplus\V^{(4)}$, hence acts
trivially on $\V^{(6)}$. This is a contradiction. Hence, we have
$\W\cap\V^{(4)}\neq\lb 0\rb$. Also, $\W\cap\V^{(4)}$ is a proper
invariant subspace of $\V^{(4)}$. Consequently, $\n^{(4)}$ is
reducible and by the case $n=4$, we have
\begin{equation}l\in\bigg\lb
r,-r^3,\unsur{r^5},\unsurr,-\unsurr\bigg\rb\end{equation} Since
$r^2\neq 1$, $(r^2)^3\neq 1$ and $(r^2)^6\neq 1$ when $\ih(6)$ is
semisimple, $(31)$ and $(32)$ imply that $l\in\lb r,-r^3\rb$. Thus,
if $\n^{(6)}$ is reducible and $l\not\in\lb\unsur{r^9},\unsur{r^3},
-\unsur{r^3}\rb$, then $l\in\lb r, -r^3\rb$. So if $\n^{(6)}$ is
reducible, then $l\in\lb r,-r^3,\unsur{r^9},\unsur{r^3},
-\unsur{r^3}\rb$.
\subsection{Proof of the main theorem when $n\geq 7$}
By the work from previous parts, the main theorem holds for $n\in\lb
3,4,5,6\rb$. When $n\geq 7$, we proceed by induction to prove the
theorem.
 Given an integer $n$ with $n\geq 7$, suppose the main theorem holds for $\n^{(n-1)}$ and $\n^{(n-2)}$. We already saw that when $l\in\lb r,-r^3,\unsur{r^{2n-3}},\unsur{r^{n-3}},-\unsur{r^{n-3}}\rb$, the representation $\n^{(n)}$ is reducible. We will show conversely that if $\n^{(n)}$ is reducible, it forces these values for $l$ and $r$. The theorem will then be proven. Suppose $\n^{(n)}$ is reducible and let $\W$ be an irreducible nontrivial proper invariant subspace of $\V^{(n)}$. By Corollary $1$, we know that $\W$ is an irreducible $\ih(n)$-module. The following proposition is part of the author's work in \cite{J}.
\begin{Proposition}
Let $K$ be a field of characteristic zero. Let $n$ be an integer with $n\geq 9$.\\
Every irreducible $K\,Sym(n)$-module is either isomorphic to one of
the Specht modules $S^{(n)}$, $S^{(n-1,1)}$, $S^{(n-2,2)}$,
$S^{(n-2,1,1)}$ or to one of their conjugates, or has dimension
greater than $\dbw$.
\end{Proposition}
\noin We have the Corollary on the dimensions:
\begin{Corollary} Assume $\ih(n)$ is semisimple. \\\vspace{-0.2in}\begin{list}{\texttt{(i)}}{}
\item Let $\D$ be an irreducible $F\,Sym(n)$-module with $n=7$ or $n\geq 9$, where $F$ is a field of characteristic zero. Then, there are two possibilities:
$$\begin{array}{lll}
\text{either} & dim\,\D\in\lb 1,n-1,\frac{n(n-3)}{2},\frac{(n-1)(n-2)}{2}\rb\\
&\\ \text{or} & dim\,\D\,>\,\frac{(n-1)(n-2)}{2}
\end{array}$$
\end{list}\vspace{-0.1in}\begin{list}{\texttt{(ii)}}{} \item \vspace{-0.1in} Let $\D$ be an irreducible $F\,Sym(8)$-module, where $F$ is a field of characteristic zero. Then $\text{dim}\,\D\in\lb 1,7,14,20,21\rb$ or $\text{dim}\,\D>21$.\end{list}
\begin{list}{\texttt{(iii)}}{}\item \vspace{-0.05in} The first two points hold if $\D$ is an irreducible $\ih(n)$-module. \end{list}
\end{Corollary}
\noin\textbf{Proof.} $(ii)$ can be seen directly by using the Hook
formula. Point $(i)$ is for $n\geq 9$ a direct consequence of
Proposition $6$ after noticing that $S^{(n-2,2)}$ has dimension
$\cil$ and $S^{(n-2,1,1)}$ dimension $\dbw$. When $n=7$, the
statement also holds by direct investigation, using for instance the
Hook formula. In characteristic zero, when the Hecke algebra
$\ih(n)$ is semisimple, the dimensions of the irreducible
$F\,Sym(n)$-modules are the same as the dimensions of the
irreducible $\ih(n)$-modules, hence $(iii)$. $\;\;\;\square$\\

\noin Let's go back to the proof of the Main Theorem. Suppose first
$n=7$ or $n\geq 9$. So $\di\,\W\in\lb 1,n-1,\cil,\dbw\rb$ or
$\di\,\W>\dbw$. First, if $\di\,\W=1$, Theorem $4$ implies that
$l=\unsur{r^{2n-3}}$. Also, if $\di\,\W=n-1$, Theorem $5$ implies
that $l\in\lb\unsur{r^{n-1}},-\unsur{r^{n-1}}\rb$. Suppose
$l\not\in\lb \unsur{r^{2n-3}},\jca, -\jca\rb$. Then we have
$\di\W\geq\cil$. 
\newtheorem{claim}{Claim}
\begin{Claim} Let $\W$ be a subspace of $\V^{(n)}$.\\
If $\di\,\W>n-1$, then $\W\cap\V^{(n-1)}\neq\lb 0\rb$.\\
If $\di\,\W>2n-3$, then $\W\cap\V^{(n-2)}\neq\lb 0\rb$
\end{Claim}
\noin\textbf{Proof.} If $\W\cap\V^{(n-1)}=\lb 0\rb$, the L-K space
$\V^{(n)}$ contains the direct sum $\W\oplus\V^{(n-1)}$, which
yields on the dimensions: $\di\,\W+\frac{(n-1)(n-2)}{2}\leq\chl$.
Then $\di\,\W\leq n-1$. Similarly, if $\W\cap\V^{(n-2)}=\lb 0\rb$,
we get\\ $\di\,\W\leq \chl-\frac{(n-2)(n-3)}{2}=2n-3$.
$\;\;\;\square$
\begin{Lemma}
When $n>6$, we have $\cil>2n-3$ and $\cil>n-1$.
\end{Lemma}
\noin By the claim and the lemma, the intersections
$\W\cap\V^{(n-1)}$ and $\W\cap\V^{(n-2)}$ are nontrivial. Since $\W$
is proper in $\V^{(n)}$, $\W$ cannot contain $\V^{(n-1)}$. Nor can
it contain $\V^{(n-2)}$. Hence $\W\cap\V^{(n-1)}$ (resp
$\W\cap\V^{(n-2)}$) is a proper nontrivial invariant subspace of
$\V^{(n-1)}$ (resp $\V^{(n-2)}$). Now $\n^{(n-1)}$ and $\n^{(n-2)}$
are both reducible. Since we assumed the main theorem to be true for
$\n^{(n-1)}$ and $\n^{(n-2)}$, we get:
$$l\in\bigg\lb r,-r^3,\unsur{r^{2n-5}},\unsur{r^{n-4}},
-\unsur{r^{n-4}}\bigg\rb\cap\bigg\lb
r,-r^3,\unsur{r^{2n-7}},\unsur{r^{n-5}},-\unsur{r^{n-5}}\bigg\rb$$
Since $r^2\neq 1$, $r^{2(n-3)}\neq 1$ and $r^{2n}\neq 1$ when
$\ih(n)$ is semisimple, it only leaves the possibility $l\in\lb
r,-r^3\rb$. \\When $n=8$, if $l\not\in\lb
\unsur{r^{13}},\unsur{r^5},-\unsur{r^5}\rb$, then we have
$\di\,\W\geq 14>13=2\times 8-3$. Hence the same method applies and
yields again $l\in\lb r,-r^3\rb$. \\Thus, we have shown that if the
representation is reducible and $l\not\in\lb
\unsur{r^{2n-3}},\jca,-\jca\rb$, then $l\in\lb r,-r^3\rb$.
$\;\;\;\square$

\section{Non semisimplicity of the BMW algebra for some specializations of its parameters}
In $\S\,2.1$, we let the Hecke algebra $H$ act on the base field $F$
by $g_i.1=r$ for all $i\in\lb 3,\dots,n-1\rb$. If we consider the
action given by $g_i.1=-\unsurr$ instead, $F$ is again a left
$H$-module for this action and we get another left $B$-module of
dimension $\chl$ by considering again the tensor product $B_1\ot F$.
We call this representation the conjugate L-K representation. By the
symmetry of the roles played by $r$ and $-\unsurr$, when $n\geq 4$
and $\ih(n)$ is semisimple, the conjugate L-K representation is
reducible exactly when $l\in\lb -\unsurr,
\unsur{r^3},-r^{2n-3},r^{n-3},-r^{n-3}\rb$. In particular, when
$n\geq 6$, since $\unsur{r^3}\not\in\lb
r,-r^3,\unsur{r^{2n-3}},\unsur{r^{n-3}}, -\unsur{r^{n-3}}\rb$, the
two representations are not equivalent. This is also true when
$n\in\lb 4,5\rb$. For instance, for the L-K representation, the
trace of the matrix of the left action by $g_{n-1}$ is
$\frac{(n-2)(n-3)}{2}\,r+\unsur{l}-(n-2)\,m$. For the conjugate
representation it is
$\frac{(n-2)(n-3)}{2}\,\big(-\unsurr\big)+\unsur{l}-(n-2)\,m$. \\We
note that Proposition $1$ remains valid for the conjugate L-K
representation. A consequence of this Proposition is that when the
representation is reducible, it is indecomposable. Then the BMW
algebra is not semisimple for the values of $l$ and $r$ for which
the L-K representation or its conjugate representation are
reducible. This is the statement of Theorem $2$.

\section{Conclusion and future developments}
In \cite{CGW}, it is established that $I_1/I_2$ is generically
semisimple where $I_1$ is the two-sided ideal $Be_1B$ and $I_2$ the
two-sided ideal generated by all the products $e_ie_j$ with
$|i-j|>2$. For each irreducible representation of the Hecke algebra
of type $A_{n-2}$ of degree $\theta$, the authors build a
generically irreducible representation of $B/I_2$ of degree
$\chl\theta$ and show that these are all the inequivalent
generically irreducible representations of $I_1/I_2$. One of the two
so-built inequivalent representations of $B/I_2$ of degree $\chl$ is
the Lawrence-Krammer representation. The other one is obtained from
the first one by replacing $r$ by $-\unsurr$. Since the
representation $B_1\ot F$ built in this paper is a generically
irreducible representation of $B/I_2$ and its kernel does not
contain $I_1$, it must be equivalent to the Lawrence-Krammer
representation. We think that the other generically irreducible
$I_1/I_2$-modules are these with $F$ replaced by an irreducible
$H$-module and that by studying these representations we could show
that $I_1/I_2$
is semisimple if and only if $l$ and $r$ don't take the specializations of the Main Theorem.\\\\




\end{document}